\documentclass[12pt,reqno]{amsart}
\usepackage{amsmath}
\usepackage{amsfonts}
\usepackage{amssymb}
\usepackage{mathtools}
\usepackage{graphicx}
\usepackage{color}
\usepackage{hyperref}
\usepackage{geometry}
\usepackage{tikz-cd}
\usepackage{enumerate}

\usepackage{xcolor}
\newtheorem{theorem}{Theorem}[section]
\newtheorem*{theorem*}{Theorem}
\newtheorem{lemma}[theorem]{Lemma}
\newtheorem{corollary}[theorem]{Corollary}
\newtheorem{proposition}[theorem]{Proposition}

\theoremstyle{definition}
\newtheorem{definition}[theorem]{Definition}
\newtheorem{example}[theorem]{Example}

\newtheorem{remark}{Remark}[section]

\DeclareMathOperator{\di}{div}

\DeclareMathOperator{\dist}{dist}
\DeclareMathOperator{\hess}{Hess}

\def\rr{{\mathbb R}}
\def\cc{{\mathbb C}}

\def\st{{S_{\mathcal{H},\widetilde{\mathcal{H}}}(f)}}
\def\dt{{df_{\mathcal{H},\widetilde{\mathcal{H}}}}}
\def\dtu{{du_{\mathcal{H},\widetilde{\mathcal{H}}}}}
\def\foli{{(M^{m+k},\mathcal{F}^k,g)}}

\def\HH{{\mathcal{H},\widetilde{\mathcal{H}}}}

\newcommand{\la}{\langle}
\newcommand{\ra}{\rangle}

\allowdisplaybreaks

\begin{document}
    \begin{abstract}
        In this paper, we consider critical points of the
        horizontal energy $E_{\HH}(f)$ for a smooth map $f$ between two Riemannian
        foliations. These critical points are referred to as horizontally harmonic
        maps. In particular, if the maps are foliated, they become transversally
        harmonic maps. By utilizing the stress-energy tensor, we establish some
        monotonicity formulas for horizontally harmonic maps from Euclidean spaces,
        the quotients $K_{m}$ of Heisenberg groups and also for transversally
        harmonic maps from Riemannian foliations with appropriate curvature pinching
        conditions. Finally, we give Jin-type theorems for either horizontally
        harmonic maps or transversally harmonic maps under some asymptotic
        conditions at infinity.
\end{abstract} %%%%%%%%%

    \title[Horizontal maps between Riemannian foliations]{On critical maps of the horizontal energy functional between Riemannian foliations}

    \author{Tian Chong}
    \date{\today}
    \address{Department of Mathematics,
Shanghai Polytechnic University, Shanghai 201209, PR China}
    \email{chongtian@sspu.edu.cn}
    
    \author{Yuxin Dong}
    \address{School of Mathematical Science  and  Laboratory of Mathematics for Nonlinear Science, Fudan University, Shanghai 200433, PR China}
    \email{yxdong@fudan.edu.cn}
    
    \author{Xin Huang }
    \address{School of Mathematics and Statistics, Nanjing University of Information Science and Technology,
Nanjing 210044,
PR China}
    \email{17110180003@fudan.edu.cn}
    
    \author{Hui Liu}
    \address{School of Mathematical Sciences, Fudan University, Shanghai 200433, PR China}
    \email{21110180015@m.fudan.edu.cn and liuh45@univie.ac.at}
    \maketitle
    
    \let\thefootnote\relax

\section{Introduction}

    Harmonic maps between Riemannian manifolds are an important object of study in
	geometric analysis. Over the decades, several generalized harmonic maps have
	appeared in the literature, playing an important role in geometry and
	topology. These generalizations include some generalized harmonic maps from
	Riemannian foliations, CR manifolds and contact manifolds (see \cite{konderakREMARKSTRANSVERSALLYHARMONIC2008, konderakTransversallyHarmonicMaps2003}, \cite{chenBochnertypeFormulasTransversally2012a}, 
	\cite{dragomirHarmonicMapsFoliated2013a},
	\cite{barlettaTransversallyHolomorphicMaps1998a}, \cite{barlettaPseudoharmonicMapsNondegenerate2001},
	\cite{petitMokSiuYeungTypeFormulas2009},
	\cite{dragomirLeviHarmonicMaps2014},
	\cite{renPseudoharmonicMapsClosed2018}, \cite{chongHARMONICPSEUDOHARMONICMAPS2019},
	etc.).
	
	Riemannian foliations are a natural generalization of Riemannian manifolds.
	Let $(M,g,\mathcal{F})\ $ and $(\widetilde{M},\widetilde{g},\widetilde{%
	\mathcal{F}})$ be two Riemannian foliations and $f:M\rightarrow \widetilde{M}
	$ be a smooth map. We say that $f$ is foliated if it maps leaves of $%
	\mathcal{F}$ to leaves of $\widetilde{\mathcal{F}}$. Then a foliated map $f$
	induces a map $\widehat{f}:M/\mathcal{F\rightarrow }\widetilde{M}/\widetilde{%
	\mathcal{F}}$ between the leaf spaces. Since the leaf spaces can be very
	complicated topological spaces, it is difficult to investigate $\widehat{f}$
	directly. One has to study related problems about the map $f$ on the upstairs
	spaces. In order to generalize the theory of harmonic maps to Riemannian
	foliations, Konderak and Wolak (cf. \cite{konderakTransversallyHarmonicMaps2003,konderakREMARKSTRANSVERSALLYHARMONIC2008}) constructed a global
	transversal tension field for a foliated map $f$ in terms of the transversal
	second fundamental form, and then they defined transversally harmonic maps
	as those smooth maps for which this tension field vanishes. However, S.
	Dragomir and A. Tommasoli (\cite{dragomirHarmonicMapsFoliated2013a}) pointed out that such maps do not, in
	general, extremize the natural transverse energy functional that will be
	described below.
	
	Let $\mathcal{V}(\mathcal{F})$ and $\mathcal{V}(\widetilde{\mathcal{F}})$ be
	the distributions that are determined by the foliations $\mathcal{F}$ and $%
	\widetilde{\mathcal{F}}$ respectively. Then we have the quotient bundles $Q(%
	\mathcal{F}):=T(M)/\mathcal{V}(\mathcal{F})$ and $Q(\widetilde{\mathcal{F}}%
	):=T(\widetilde{M})/\mathcal{V}(\widetilde{\mathcal{F}})$, which are called
	the transverse bundles of the foliations. These transverse bundles get
	natural fiberwise metrics from $g$ and $\widetilde{g}$ respectively. The
	differential map $df$ of $f$ induces the transverse differential $d_{T}f:Q(%
	\mathcal{F})\rightarrow f^{-1}(Q(\widetilde{\mathcal{F}}))$ in a natural
	way. Consequently we may define the norm $\Vert d_{T}f\Vert $ and the
	following transverse energy functional $E_{T}(f)$ (cf. \cite{dragomirHarmonicMapsFoliated2013a} for details): 
	\begin{equation}\label{tag1}
	E_{T}(f)=\int_{M}e_{T}(f)dV_{g}.  
	\end{equation}%
	It turns out that the transversally harmonic maps defined by Konderak and
	Wolak are not the critical points of the functional $E_{T}(f)$, but rather
	the critical points of the following energy functional%
	\begin{equation}\label{tag2}
	E_{T}^{\ast }(f)=\int_{M}\frac{e_{T}(f)}{vol_{L}}dV_{g} , \tag{2}
	\end{equation}%
	where $vol_{L}$ denotes the volume of each fiber. In \cite{dragomirHarmonicMapsFoliated2013a}, S. Dragomir and
	A. Tommasoli proposed a different definition of generalized harmonic maps,
	which are exactly the extremals of $E_{T}(f)$ through any variation of $f$
	by foliated maps. They call such maps $(\mathcal{F},\widetilde{\mathcal{F}})$%
	-harmonic maps. It is easy to verify that the above two definitions of
	generalized harmonic maps coincide when the foliation $\mathcal{F}$ is
	harmonic (cf. \cite{chenBochnertypeFormulasTransversally2012a}, \cite{dragomirHarmonicMapsFoliated2013a}).
	
	In this paper, we introduce the horizontal energy functional $E_{\mathcal{H},%
	\widetilde{\mathcal{H}}}(f)$ similar to \eqref{tag1} for any map $f:(M^{m+k},\mathcal{%
	F},g)\rightarrow (N^{n+l},\widetilde{\mathcal{F}},\widetilde{g})$ (not
	necessarily foliated) between the Riemannian foliations, and then define
	horizontally harmonic maps as the critical points of $E_{\mathcal{H},%
	\widetilde{\mathcal{H}}}(f)$ through any variation (see \S \ref{se:criticalmaps} for details).
	The notion of horizontally harmonic maps slightly generalizes the notion of $%
	(\mathcal{F},\widetilde{\mathcal{F}})$-harmonic maps in \cite{dragomirHarmonicMapsFoliated2013a}. Notice that
	removing the foliated condition about the map in our definition allows the
	notion to involve the subelliptic harmonic maps when the domain is a
	Riemannian foliation and the target is a Riemannian manifold with the
	trivial foliation by points. The main purpose of the present paper is to
	study the properties of the horizontally harmonic maps, especially the
	Liouville type theorems for them.
	
	As known, the stress-energy tensors are a useful tool for deriving
	monotonicity formulas of harmonic maps or their generalizations (cf. \cite{price1983amono},
	\cite{xin1986diff}, \cite{dongVanishingTheoremsVector2011a}, \cite{dlyLiou2016}, and the references therein). Since the stress-energy tensors
	are naturally linked to the conservation laws of related energy functionals,
	the monotonicity formulas of the energies follow from Stokes' theorem,
	coarea formula and comparison theorems in Riemannian geometry. These
	monotonicity formulas can be used to establish Liouville type theorems in
	the following two ways. One way is to derive Liouville type theorems
	directly from the monotonicity formulas by assuming the energy growth
	conditions (cf. \cite{sealeysome1982}, \cite{hu1984anonexist}, \cite{ll2004theenergy}, \cite{dongVanishingTheoremsVector2011a}, \cite{rs2000energy}, and the references therein).
	Another way is to establish Liouville type theorems by assuming suitable
	asymptotic conditions of the maps at infinity (cf. \cite{jin_liouville_nodate}, \cite{dlyLiou2016}, \cite{rs2000energy}).
	
	For our purpose, we introduce the stress-energy tensor $S_{\mathcal{H},%
	\widetilde{\mathcal{H}}}(f)$ associated with the energy functional $E_{%
	\mathcal{H},\widetilde{\mathcal{H}}}(f)$ and derive the corresponding
	conservation law formula for $(\di S_{\mathcal{H},\widetilde{\mathcal{H
    	}}}(f))(X)$ with respect to any vector field $X\in \mathfrak{X}(M)$. It
	turns out that a horizontally harmonic map $f:(M^{m+k},\mathcal{F}
	,g)\rightarrow (N^{n+l},\widetilde{\mathcal{F}},\widetilde{g})$ does not
	satisfy the conservation law in general, but if $f$ is a transversally
	harmonic map, then it satisfies the conservation law. Even though in the
	general case $f$ does not satisfy the conservation law, we find that the
	conservation law formula becomes more manageable if  the domain
	Riemannian foliation is simple. To illustrate this, we consider the
	following two cases for horizontally harmonic maps: (i) $(M^{m+k},\mathcal{F%
	},g)=(%
	\mathbb{R}
	^{m+k},%
	\mathbb{R}
	^{k},g)$ the Euclidean space foliated by Euclidean subspaces, where $g$ is
	some mixed conformally flat metric on $%
	\mathbb{R}
	^{m+k}$ (see Example \ref{exam4}); (ii) $(M^{m+k},\mathcal{F},g)=(K_{m},S^{1},g)$
	the quotient space of the Heisenberg group $H_{m}$ with the canonical metric 
	$g$ (see Example \ref{ex:sasa}). The horizontal distribution of the former is
	integrable, while the horizontal distribution of the latter is
	non-integrable. They can be viewed respectively as the simplest models in
	the two cases. Since the conservation law formula becomes simpler in
	horizontal directions for these cases, it is more appropriate to apply the
	formula on cylindrical regions defined by the distance function from a fixed
	leaf. Consequently, we can establish some monotonicity formulas of the
	horizontal energy on cylindrical regions for horizontally harmonic maps from
	either $\rr^{m+k}$ or $K_{m}$. As mentioned above, a transversally harmonic
	map  $f:(M^{m+k},\mathcal{F},g)\rightarrow (N^{n+l},\widetilde{\mathcal{F}},%
	\widetilde{g})$ between two Riemannian foliations satisfies the conservation
	law, that is, $(\di S_{\mathcal{H},\widetilde{\mathcal{H}}}(f))(X)=0$
	in any direction $X$. This enables us to establish some monotonicity
	formulas of the horizontal energy on geodesic balls for transversally
	harmonic maps under suitable curvature pinching conditions of $(M,g)$.
	Obviously, under the assumption of appropriate energy growth conditions, the
	above monotonicity formulas for horizontally harmonic maps or transversally
	harmonic maps give us directly Liouville type theorems for these maps.
	
	In \cite{jin_liouville_nodate}, Jin established several interesting Liouville type theorems for
	harmonic maps from Euclidean spaces endowed with conformally flat metrics,
	under an asymptotic condition of these maps at infinity. A special case of
	his results is that if $u:(\rr^{m},g_{can})\rightarrow (N^{n},h)$ is a
	harmonic map, and $u(x)$ converges to a fixed point $p_{0}\in N$ as $|
	x| \rightarrow \infty $, then $u$ is a constant map. His method is based
	on the following two steps for a non-constant harmonic map: first, deriving
	a lower bound on the energy growth rate, and second, deriving an upper bound
	on the energy growth rate under an asymptotic condition. If the two growth
	rates are incompatible, this shows that $u$ can only be constant. Inspired
	by Jin's method, we next investigate Jin-type Liouville theorems for a
	horizontally harmonic map or transversally harmonic map $f:(M^{m+k},\mathcal{%
	F},g)\rightarrow (N^{n+l},\widetilde{\mathcal{F}},\widetilde{g})$ by
	assuming that the map approaches a fixed point $q$ at infinity. We will
	establish Jin-type theorems for the following two cases: (i) $f:(
	\mathbb{R}
	^{m+k},
	\mathbb{R}
	^{k},g)\rightarrow (N^{n+l},\widetilde{\mathcal{F}},\widetilde{g})$ is a
	horizontally harmonic map, where $g$ is a mixed conformally flat metric on $
	\mathbb{R}
	^{m+k}$; (ii) $f:(M^{m+k},\mathcal{F},g)\rightarrow (N^{n+l},\widetilde{%
	\mathcal{F}},\widetilde{g})$ is a transversally harmonic map, where $g$ is a
	complete Riemannian metric with radial curvature $K_{r}$ satisfying $-\alpha
	^{2}\leq K_{r}\leq -\beta ^{2}$, $\alpha ,\beta >0$ and $(m+k-1)\beta
	-2\alpha \geq 0$. Note that the monotonicity formulas we have already
	established for these two cases in \S 5 actually give lower bounds on the
	energy growth rate. This completes the first step in Jin's method. For the
	second step, we may assume that the image of the map $f$ near infinity is
	contained in a foliated coordinate chart $(U,\widetilde{\mathcal{F}}|_{U})$,
	which induces a natural Riemannian submersion $\pi _{U}:U\rightarrow 
	\mathcal{B}_{U}$. Clearly, estimating the horizontal energy of $f$ near
	infinity is equivalent to estimating the horizontal energy of $\pi _{U}\circ
	f$. Then we can apply Jin's method to give the upper bound on the energy
	growth rate. Therefore, if the two growth rates are incompatible, we can
	establish a Jin-type Liouville theorem (see \S \ref{se:Jintypethm} for details).

This paper is structured in the following manner. 
In \S \ref{se:Riemfoli}, we collect some basic notions, propositions and some examples in Riemannian foliations. 
In \S \ref{se:criticalmaps}, we introduce the horizontal energy $E_{\mathcal{H}, \widetilde{\mathcal{H}}}(f)$ for a smooth map $f$, and then deduce the first variation formula, and we also explore the geometric meaning of horizontally constant maps.
In \S \ref{se:streeenergy}, we derive the divergence formula for the stress-energy tensor $S_{\HH}(f)$ and introduce the associated concepts of conservation laws.   
In \S \ref{se:monoformula}, 
we establish monotonicity formulas for horizontally harmonic maps from Euclidean spaces and the quotient space $K_m$ of the Heisenberg group $H_m$, as well as for transversally harmonic maps from Riemannian foliations under certain curvature pinching conditions.
 Finally, in \S \ref{se:Jintypethm}, we present Jin-type theorems for horizontally harmonic maps from mixed conformally flat Euclidean spaces and transversally harmonic maps from  Riemannian foliations with a pole.

~

\section{Riemannian foliations}\label{se:Riemfoli}

We first recall some notions and notations in foliation theory. Let $\mathcal{V}$ be an integrable $k$-dimensional distribution on a smooth manifold $M^{m+k}$ of dimension $m+k$. The collection $\mathcal{F}^k$ of integral submanifolds of $\mathcal{V}$ is called a $k$-dimensional foliation or a foliation of codimension $m$ on $M$. The pair $\left(M^{m+k}, \mathcal{F}^k\right)$ is said to be a foliated manifold. Clearly $\mathcal{F}^k$ gives a partition of $M$ into disjoint $k$ dimensional immersed submanifolds $\mathcal{L}_\alpha$, which are called the leaves of $\mathcal{F}^k$. The space of leaves, denoted by $\mathcal{B}$, is the quotient space of the equivalence relation: $x \sim y$ if $x$ and $y$ lie on the same leaf. Notice that, in the case of an arbitrary foliation, $\mathcal{B}$ is possibly not Hausdorff. By Frobenius theorem, every point of $M$ has a foliated coordinate chart $(U, \varphi)$ with coordinates $\left(x_1, \ldots, x_m, y_1, \ldots, y_k\right)$ such that $\operatorname{span}\left\{\frac{\partial}{\partial y_\alpha}\right\}_{1 \leq \alpha  \leq k}=\left.\mathcal{V}\right|_U$. Let $D_\rho^{m+k}$ denote the open neighborhood of the origin in $\mathbb{R}^{m+k}$ given by
$$
\begin{aligned}
    D_\rho^{m+k} & =\left\{\left(x_1, \ldots, x_m, y_1, \ldots, y_k\right) \in \mathbb{R}^{m+k} ~|~ | x_i|<\rho,| y_\alpha |<\rho\right\} \\
    & =D_\rho^m \times D_\rho^k
\end{aligned}
$$
for some $\rho>0$. Without loss of generality, one may assume that $\varphi(U)=D_\rho^{m+k}$. Clearly $\left.\mathcal{F}^k\right|_U$ induces a local submersion $\pi_U: U \rightarrow \mathcal{B}_U=\left(\left.\mathcal{F}^k\right|_U\right) / \sim$, which corresponds to the natural projection
\begin{equation}
    \begin{aligned}
        \operatorname{Pr}_m: D_\rho^{m+k} & \rightarrow D_\rho^m \\
        \left(x_1, \ldots, x_m, y_1, \ldots, y_k\right) & \longmapsto\left(x_1, \ldots, x_m\right) .
    \end{aligned}
\end{equation}
If $\mathcal{F}^k$ induces a global submersion from $M$ to a smooth base manifold $\mathcal{B}$, it is often called a simple foliation. When the foliation consists of 0-dimensional points, $(M, \mathcal{F})$ coincides with $M$, in which case we say the foliation $(M, \mathcal{F})$ is a point foliation.

We are interested in transverse structures of foliations. Let $\left(M^{m+k}, \mathcal{F}^k\right)$ be a foliated manifold with $k$-dimensional integrable subbundle $\mathcal{V}$. There is an exact sequence of vector bundles
\begin{equation}
    0 \longrightarrow \mathcal{V} \longrightarrow T M \xrightarrow{\pi_Q} Q \longrightarrow 0
\end{equation}
where the quotient bundle $Q=T M / \mathcal{V}$ is called the normal bundle of the foliation. A vector field $X$ on $M$ is said to be foliated with respect to $\mathcal{F}^k$ if $[V, X]$ is tangent to the leaves of $\mathcal{F}^k$ for any vector field $V$ tangent to the leaves of $\mathcal{F}$. Equivalently, the local 1-parameter group of $X$ preserves the foliation. It is easy to show that in a foliated coordinate chart $\left(U ; x_1, \ldots, x_m, y_1, \ldots, y_k\right)$, a foliated vector field $X$ can be expressed as
\begin{equation}
    X=\sum_{\alpha=1}^k A^{\alpha }\left(x_1, \ldots, x_m, y_1, \ldots, y_k\right) \frac{\partial}{\partial y_\alpha}+\sum_{i=1}^m B^{i}\left(x_1, \ldots, x_m\right) \frac{\partial}{\partial x_i},
\end{equation}
where $A^\alpha$ and $B^i$ are smooth functions of the listed variables. Clearly $\left.X\right|_U$ projects to a vector field on $\mathcal{B}_U$, given by
$$
\bar{X}=\sum_{i=1}^m B^{i}\left(x_1, \ldots, x_m\right) \frac{\partial}{\partial x_i}.
$$
Hence foliated vector fields are locally projectable vector fields.

We now endow the foliated manifold $\left(M^{m+k}, \mathcal{F}^k\right)$ with a Riemannian metric
$g$. Then the tangent bundle of $M$ admits an orthogonal decomposition
\begin{equation}\label{eq01}
    T M=\mathcal{H} \oplus \mathcal{V},
\end{equation}
where $\mathcal{H}=\mathcal{V}^{\perp}$ is the orthogonal complement of $\mathcal{V}$ with respect to $g$. The subbundles $\mathcal{H}$ and $\mathcal{V}$ are called the horizontal and vertical bundles of $\mathcal{F}$, respectively. Let $\mathfrak{X}(M)$ represent the space of smooth vector fields on $M$. The sets of smooth sections of $\mathcal{H}$ and $\mathcal{V}$ are denoted by $\mathfrak{X}_{\mathcal{H}}$ and $\mathfrak{X}_{\mathcal{V}}$, respectively. An element in $\mathfrak{X}_{\mathcal{H}}$ (respectively $\mathfrak{X}_{\mathcal{V}}$ ) is known as a horizontal (respectively vertical) vector field. Obviously, restricting $\pi_Q$ to $\mathcal{H}$ gives a bundle isomorphism $\tau: \mathcal{H} \rightarrow Q$. According to the decomposition $\eqref{eq01}$, we also have the following natural projections
\begin{equation}
    \pi_{\mathcal{H}}: T M \rightarrow \mathcal{H}, \quad \pi_{\mathcal{V}}: T M \rightarrow \mathcal{V},
\end{equation}
which are called the horizontal and vertical projections respectively. Defining
\begin{equation}
    g_{\mathcal{H}}=g\left(\pi_{\mathcal{H}}, \pi_{\mathcal{H}}\right), \quad g_{\mathcal{V}}=g\left(\pi_{\mathcal{V}}, \pi_{\mathcal{V}}\right),
\end{equation}
one may express $g$ as
\begin{equation}
    g=g_{\mathcal{H}}+g_{\mathcal{V}} .
\end{equation}
Following \cite{reinhartFoliatedManifoldsBundleLike1959}, a Riemannian metric $g$ on $(M, \mathcal{F})$ is said to be bundle-like if the Lie derivative $\mathcal{L}_V g_{\mathcal{H}}=0$ for any vertical vector field $V$. Equivalently, on each foliated coordinate chart $U$ with the local submersion $\pi_U: U \rightarrow \mathcal{B}_U$, there exists a Riemannian metric $h$ on $\mathcal{B}_U$ for which the local submersion $\pi_U:(U, g) \rightarrow$ $\left(\mathcal{B}_U, h\right)$ becomes a Riemannian submersion. In this case, we say that $(M, \mathcal{F}, g)$ is a Riemannian foliation. Traditionally the terminology  ``Riemannian foliation" refers only to a foliation endowed with a transverse Riemannian metric on the quotient bundle. It is a known fact (cf. \cite{molinoRiemannianFoliations1988}) that, if $g$ is a bundle-like metric on $(M, \mathcal{F})$, then it induces an associated transverse Riemannian metric $g_Q$ on $Q$ from $g_{\mathcal{H}}$ through the isomorphism $\tau: \mathcal{H} \rightarrow Q$. Conversely, if the quotient bundle $Q$ has a transverse Riemannian metric $g_Q$, then there is also a bundle-like metric $g$ that produces the same transverse metric on $Q$. Thus, in this paper, the triple $(M, \mathcal{F}, g)$ is referred to as a Riemannian foliation when $g$ is a bundle-like metric. In particular, if $\mathcal{F}$ is a simple foliation, it corresponds to a global Riemannian submersion from $(M, g)$ to a base manifold $\mathcal{B}$.

Let $\left(M^{m+k}, \mathcal{F}^k , g\right)$ be a Riemannian foliation and let $\nabla$ be the Levi-Civita connection of $g$. Following \cite{gromollMetricFoliationsCurvature2009a}, we introduce two $\mathcal{V}$-valued fundamental tensor fields $\mathcal{W}$ and $\mathcal{T}$ on $M$ that measure the complexity of the Riemannian foliation (notice that our notations differ from those in \cite{gromollMetricFoliationsCurvature2009a}). The $\mathcal{W}$-tensor is the tensor field $\mathcal{W}: \mathcal{H} \times \mathcal{V} \rightarrow \mathcal{V}$ defined by
\begin{equation}\label{wten}
    \mathcal{W}_X U=-\pi_{\mathcal{V}}\left(\nabla_U X\right), \quad X \in \mathfrak{X}_{\mathcal{H}}, \quad U \in \mathfrak{X}_{\mathcal{V}},
\end{equation}
and the $\mathcal{T}$-tensor is the tensor field $\mathcal{T}: \mathcal{H} \times \mathcal{H} \rightarrow \mathcal{V}$ given by
\begin{equation}\label{ttensor}
    \mathcal{T}_X Y=\pi_{\mathcal{V}}\left(\nabla_X Y\right)=\frac{1}{2} \pi_{\mathcal{V}}([X, Y]), X, Y \in \mathfrak{X}_{\mathcal{H}} .
\end{equation}		
The second equality in $\eqref{ttensor}$ is shown in \cite{gromollMetricFoliationsCurvature2009a} by using the Koszul formula (see the proof of Theorem 1.2.1 in \cite{gromollMetricFoliationsCurvature2009a}). This implies that $\mathcal{T}$ is anti-symmetric with respect to $X$ and $Y$. Notice that $\mathcal{W}_X$ is just the Weingarten transformation of a leaf in the direction $X$. If we define an $\mathcal{H}$-valued tensor $\mathcal{S}: \mathcal{V} \times \mathcal{V} \rightarrow \mathcal{H}$ by
\begin{equation}\label{eq10}
    \mathcal{S}(U, V)=\pi_{\mathcal{H}}\left(\nabla_U V\right), \quad U, V \in \mathfrak{X}_{\mathcal{V}}
\end{equation}
then
\begin{equation}\label{eq11}
    g\left(\mathcal{W}_X U, V\right)=g(\mathcal{S}(U, V), X) .
\end{equation}		
Consequently $\mathcal{W} \equiv 0$, or equivalently $\mathcal{S} \equiv 0$ if and only if the leaves are totally geodesic, in which case the Riemannian foliation is referred to as totally geodesic. Clearly $\mathcal{T} \equiv 0$ if and only if the horizontal distribution $\mathcal{H}$ is integrable.

One may extend $\mathcal{W}, \mathcal{S}$ and $\mathcal{T}$ to the entire tangent bundle $T M$ by defining
\begin{equation}\label{eq12}
    \begin{aligned}
        \mathcal{W}_X Y & =\mathcal{W}_{\pi_{\mathcal{H}}(X)} \pi_{\mathcal{V}}(Y), \quad \mathcal{S}(X, Y)=\mathcal{S}\left(\pi_{\mathcal{V}}(X), \pi_{\mathcal{V}}(Y)\right) \\
        \mathcal{T}_X Y & =\mathcal{T}_{\pi_{\mathcal{H}}(X)} \pi_{\mathcal{H}}(Y)
    \end{aligned}
\end{equation}
for any $X, Y \in T M$. Then $\mathcal{W}, \mathcal{S}$ and $\mathcal{T}$ become smooth tensor fields on $M$. The mean curvature vector field of the foliation is a global horizontal vector field $\kappa$ defined by
\begin{equation}
\kappa=\operatorname{tr}_g \mathcal{S}.
\end{equation}
According to \cite{gromollMetricFoliationsCurvature2009a}, the tensor fields $\mathcal{W}$ and $\mathcal{T}$ essentially determine the geometry of the Riemannian foliation. From the geometry of submanifolds, we know that $\mathcal{W}_X$ is self-adjoint with respect to $g$ for any $X \in T M$. The adjoint of $\mathcal{T}_X$ for any $X \in T M$ is defined by (\cite{gromollMetricFoliationsCurvature2009a})
\begin{equation}\label{eq14}
    g\left(\mathcal{T}_X^* Y, Z\right)=g\left(Y, \mathcal{T}_X Z\right), \quad Y, Z \in T M .
\end{equation}
Using $\eqref{ttensor}$, $\eqref{eq12}$ and $\eqref{eq14}$, it is easy to see that $\mathcal{T}^*$ is an $\mathcal{H}$-valued tensor field on $M$.

For local computations, it is convenient to choose a local orthonormal frame field $\left\{e_A\right\}_{A=1}^{m+k}$ in $\left(M^{m+k}, \mathcal{F}^k , g\right)$ such that
$$
\operatorname{span}\left\{e_i\right\}_{i=1}^m=\mathcal{H}, \quad \operatorname{span}\left\{e_\alpha\right\}_{\alpha=m+1}^{m+k}=\mathcal{V} .
$$
Such a frame field is referred to as an adapted frame field for the foliation. Using an adapted frame field, the mean curvature vector can be expressed by
\begin{equation}
    \kappa=\sum_{\alpha=m+1}^{m+k} \mathcal{S}\left(e_\alpha, e_\alpha\right)=\sum_{\alpha=m+1}^{m+k} \pi_{\mathcal{H}}\left(\nabla_{e_\alpha} e_\alpha\right).
\end{equation}
Recall that a basic vector field in $(M, \mathcal{F} , g)$ is one that is both horizontal and foliated. From the above discussion, we know that there is a local Riemannian submersion $\pi_U:(U, g) \rightarrow\left(\mathcal{B}_U, h\right)$ around each point $p$ in $M$. If $\bar{X}$ is a vector field on $\mathcal{B}_U$, then there exists a unique basic lift of $\bar{X}$, which is a smooth vector field $X$ on $U$ such that $X$ is horizontal and $d \pi_U(X)=\bar{X}$. Let $\left\{\bar{e}_i\right\}_{i=1}^m$ be a local orthonormal frame field around $\bar{p}=\pi_U(p)$ in $\left(\mathcal{B}_U, b\right)$, and let $\left\{e_i\right\}_{i=1}^m$ be the basic lift of $\left\{\bar{e}_i\right\}_{i=1}^m$ with respect to $\pi_U$. Then $\left\{e_i\right\}_{i=1}^m$ is a local orthonormal frame field of $\mathcal{H}$ in $U$.

\begin{lemma}[\cite{gromollMetricFoliationsCurvature2009a}]\label{lem1}
    Let $(M, \mathcal{F} , g)$ be a Riemannian foliation with Levi-Civita connection $\nabla$. If $X, Y \in \mathfrak{X}(M)$ are basic, then so is $\pi_{\mathcal{H}}\left(\nabla_X Y\right)$. In fact, if $\bar{X}=d \pi_U(X)$ and $\bar{Y}=d \pi_U(Y)$ for the local Riemannian submersion $\pi_U$ : $(U, g) \rightarrow\left(\mathcal{B}_U, h\right)$, then $d \pi_U\left(\nabla_X Y\right)=(\nabla^{\mathcal{B}}_{\bar{X}} \bar{Y}) \circ \pi$, where $\nabla^{\mathcal{B}}$ denotes the Levi-Civita connection of $\left(\mathcal{B}_U, h\right)$.
\end{lemma}

One can always find an orthonormal frame field $\left\{\bar{e}_i\right\}_{i=1}^m$ in the Riemannian manifold $\left(\mathcal{B}_U, h\right)$ such that $\left(\nabla^{\mathcal{B}} \bar{e}_i\right)_{\bar{p}}=0$. It follows from Lemma \ref{lem1} that around any point $p \in M$, there exists an adapted frame field $\left\{e_A\right\}_{A=1}^{m+k}$ such that $e_i\, (1 \leq$ $i \leq m)$ is basic, and
\begin{equation}
    \pi_{\mathcal{H}}\left(\nabla_X e_i\right)_p=0, \quad 1 \leq i \leq m,
\end{equation}
for any $X \in \mathcal{H}_p$. Such an adapted frame field will be useful for computations in Riemannian foliations. We shall also need the following lemma. For the convenience of the readers, we present its simple proof given by \cite{gromollMetricFoliationsCurvature2009a}.

\begin{lemma}[\cite{gromollMetricFoliationsCurvature2009a}]\label{le:td}
    Let $(M, \mathcal{F} , g)$ be a Riemannian foliation with Levi-Civita connection $\nabla$. If $X \in \mathfrak{X}_{\mathcal{H}}$ is basic and $Y \in \mathfrak{X}_{\mathcal{V}}$, then $\mathcal{T}_X^* Y=-\pi_{\mathcal{H}}\left(\nabla_Y X\right)$.
\end{lemma}
\begin{proof}
    Since $X$ is basic and $Y \in \mathfrak{X}_{\mathcal{V}}$, we have $[X, Y] \in \mathfrak{X}_{\mathcal{V}}$. For any $Z \in \mathfrak{X}_{\mathcal{H}}$, we deduce by using $\eqref{ttensor}$ and $\eqref{eq14}$ that
    $$
    \begin{aligned}
        \left\langle\mathcal{T}_X^* Y, Z\right\rangle & =\left\langle Y, \mathcal{T}_X Z\right\rangle=\left\langle Y, \nabla_X Z\right\rangle \\
        & =X\langle Y, Z\rangle-\left\langle\nabla_X Y, Z\right\rangle \\
        & =-\left\langle\nabla_Y X, Z\right\rangle-\langle[X, Y], Z\rangle \\
        & =-\left\langle\nabla_Y X, Z\right\rangle .
    \end{aligned}
    $$
    This proves the lemma.
\end{proof}

Before concluding this section, we would like to give some examples of Riemannian foliations.

\begin{example}\label{exam3}
    Let $\mathbb{R}^{m+k}$ be the $(m+k)$-dimensional Euclidean space with the canonical Euclidean metric $g_{c a n}$ and $\pi: \mathbb{R}^{m+k} \rightarrow \mathbb{R}^m$ be the natural projection. Set $\mathcal{F}=\{\pi^{-1}(q)=\mathbb{R}^k: \forall q \in \mathbb{R}^m\}$. Then $(\mathbb{R}^{m+k}, \mathcal{F} , g_{c a n})$ is a Riemannian foliation with $\mathcal{W} \equiv 0$ and $\mathcal{T} \equiv 0$.
\end{example}

\begin{example}\label{exam4}
Let us consider a mixed conformal transformation of
the Euclidean metric $g_{can}$ in Example \ref{exam3} as follows:
\[
g=\phi (x)g_{can}^{h}+\eta (x,y)g_{can}^{v}
\]
with $\phi ,\eta >0$, where $g_{can}^{h}$ and $g_{can}^{v}$ are the
canonical metrics of $
\mathbb{R}
^{m}$ and $
\mathbb{R}
^{k}$ respectively. It is easy to see that the natural projection $\pi :(
\mathbb{R}
^{m+k},
\mathbb{R}
^{k},g)\rightarrow 
\mathbb{R}
^{m}$ is a simple Riemannian foliation.	
\end{example}

Let $\mathcal{H}$ be a distribution on $M$. A Lipschitz curve $\gamma:[0, l] \rightarrow M$ is called horizontal if $\gamma^{\prime}(t) \in \mathcal{H}_{\gamma(t)}$ a.e. in $[0, l]$. We say that $\mathcal{H}$ satisfies the Hörmander's condition of order $r$, if sections of $\mathcal{H}$ together with their Lie brackets up to order $r$ span $T_x M$ at each point $x$. If $\mathcal{H}$ satisfies the Hörmander's condition, we know from the theorem of Chow-Rashevsky (\cite{chowUeberSystemeLiearren1940},  \cite{rashevskii1938two}) that there always exist horizontal curves joining any two points $p_1$ and $p_2$ in $M$.

\begin{example}\label{ex:sasa}
    From \cite{blairRiemannianGeometryContact2010a} and   \cite{dragomirDifferentialGeometryAnalysis2006}, we know that any Sasakian
    manifold admits a Riemannian foliation, whose leaves are integral curves of
    the Reeb vector field. The horizontal distribution of any Sasakian manifold
    satisfies the H\"{o}rmander's condition. The odd dimensional sphere $
    (S^{2m+1},g)\subset (
    \mathbb{C}
    ^{m+1},g_{can})\ $with induced CR structure
     and standard metric is one of the
    simplest Sasakian manifolds. 
    Another important example is the Heisenberg group $H_{m}=
\mathbb{C}^{m}\times \rr$. It is a Lie group with the following group law:
\[
(z,t)\cdot (w,s)=(z+w,t+s+2\operatorname{Im}(\sum_{j=1}z^{j}\overline{w}^{j})),
\]%
where $(z,t)=(z^{1},...,z^{m},t)$, $(w,s)=(w^{1},...,w^{m},s)\in 
\mathbb{C}
^{m}\times \rr$. Write $z^{j}=x^{j}+\sqrt{-1}y^{j}$ for $j=1,...,m$. Then $%
(x^{1},...,x^{m},y^{1},...,y^{m},t)$ is a real global coordinate system of $%
H_{m}$. Set%
\[
\eta =\frac{1}{2}dt+\sum_{j=1}^{m}\left( x^{j}dy^{j}-y^{j}dx^{j}\right) .
\]%
The canonical Riemannian metric (Webster metric) on $H_{m}$ is given by 
\[
g=\sum_{j=1}^{m}\left( (dx^{j})^{2}+(dy^{j})^{2}\right) +\eta \otimes \eta .
\]%
The natural projection $\pi :H_{m}=%
\mathbb{C}
^{m}\times \rr\rightarrow 
\mathbb{C}
^{m}=%
\mathbb{R}
^{2m}$ is a Riemannian submersion with totally geodesic leaves, where $%
\mathbb{R}
^{2m}$ is endowed with the canonical Euclidean metric $g_{can}$. The
vertical vector field $\xi =\partial /\partial t$ is a Killing vector field
on $H_{m}$ whose $1$-parameter transformation group $\{\exp (t\xi )\}_{t\in 
\mathbb{R}
}$ preserves the foliation (cf. \cite{dragomirDifferentialGeometryAnalysis2006},  \cite{gromollMetricFoliationsCurvature2009a}). Set $K_{m}=H_{m}\diagup
\{\exp (k \xi)\}_{k\in 
\mathbb{Z}
}=
\mathbb{C}
^{m}\times S^{1}$. Then the induced projection $\pi :K_{m}\rightarrow 
\mathbb{C}
^{m}$ is also a Riemannian submersion with totally geodesic leaves.
\end{example}

\begin{example}\label{ex:tang}
     Let $(M^m,g)$ be a Riemannian manifold and $\nabla $ be
    its Levi-Civita connection. Let $\pi :TM\rightarrow M$ denote the natural
    projection from the tangent bundle $TM$ to $M$. The kernel of $d\pi $ gives
    a vertical distribution $\mathcal{V}$ on $TM$, that is, $\mathcal{V}%
    _{(p,u)}=\ker (d\pi | _{(p,u)})$ for any $(p,u)\in TM$. In terms of the
    connection $\nabla $, a horizontal distribution $\mathcal{H}$ on $TM$ can be
    defined in the usual way. As a result, we have a direct sum%
    \[
    T_{(p,u)}(TM)=\mathcal{H}_{(p,u)}\oplus \mathcal{V}_{(p,u)} 
    \]%
    for any $(p,u)\in TM$. It is known that $g$ determines a natural Riemannian
    metric $\widehat{g}$ on $TM$, called Sasaki metric, such that $\mathcal{H}=%
    \mathcal{V}^{\bot }$ with respect to $\widehat{g}$ and $\pi :(TM,\widehat{g}%
    )\rightarrow (M,g)$ is a Riemannian submersion with totally geodesic leaves $%
    \{\pi ^{-1}(p)=T_{p}M$, $p\in M\}$ (cf. \cite{sasakiDifferentialGeometryTangent1958}, \cite{yanoTangentCotangentBundles1973},  \cite{gudmundssonGeometryTangentBundles2002}).
Next, let $T_{1}M$ denote the unit tangent sphere
bundle, consisting of the unit tangent vectors on $(M,g)$. Then $T_{1}M$ has
a Riemannian metric $g_{S}$ induced from the Sasaki metric $\widehat{g}$
(cf. \cite{Boeckx1997CharacteristicRO,boeckxUnitTangentSphere2001}). The natural projection $\pi
_{1}:(T_{1}M,g_{S})\rightarrow (M,g)$ is also a Riemannian submersion with
totally geodesic leaves $\{\pi _{1}^{-1}(p)=S^{m-1}(T_{p}M)$, $p\in M\}$. In
terms of the results on the Lie bracket of horizontal vector fields on $TM$
(cf. \cite{dombrowskiGeometryTangentBundle1962}, \cite{gudmundssonGeometryTangentBundles2002}) and the Lie bracket of horizontal vector fields on $%
T_{1}M$ (cf. \cite{boeckxUnitTangentSphere2001}), it is possible to construct examples of tangent
bundles and unit tangent bundles whose horizontal distributions satisfy the H%
\"{o}rmander's condition.
\end{example}
~

\section{Critical maps of the horizontal energy functional}\label{se:criticalmaps}

Let $(M^{m+k}, \mathcal{F}^k , g)$ and $(N^{n+l}, \widetilde{\mathcal{F}}^l , \widetilde{g})$ be two Riemannian foliations. We shall follow the notations in the previous section for Riemannian foliations, and denote the corresponding geometric data of $N$, such as the horizontal and vertical distributions, the Levi-Civita connection, the fundamental tensors, etc., by the same notations as in $M$, but with $\sim$ on them. For simplicity, we shall often use $\langle\cdot, \cdot\rangle$ to denote the inner products induced from $g$ or $\widetilde{g}$ on various tensor bundles on $M$ or $N$.

Let $f: M \rightarrow N$ be a smooth map between $M$ and $N$. Notice that the identity isomorphisms of the tangent bundles $T M$ and $T N$ can be expressed as
$$
i d_{T M}=\pi_{\mathcal{H}}+\pi_{\mathcal{V}}, \quad i d_{T N}=\pi_{\widetilde{\mathcal{H}}}+\pi_{\widetilde{\mathcal{V}}},
$$
respectively. Hence the differential $d f$ decomposes into the following partial differentials
\begin{equation}\label{dfde}
    \begin{aligned}
        d f & =\left(\pi_{\widetilde{\mathcal{H}}}+\pi_{\widetilde{\mathcal{V}}}\right) \circ d f \circ\left(\pi_{\mathcal{H}}+\pi_{\mathcal{V}}\right) \\
        & =d f_{\mathcal{H}, \widetilde{\mathcal{H}}}+d f_{\mathcal{V}, \widetilde{\mathcal{H}}}+d f_{\mathcal{H}, \widetilde{\mathcal{V}}}+d f_{\mathcal{V}, \widetilde{\mathcal{V}}} \\
        & =d f_{\mathcal{H}}+d f_{\mathcal{V}} \\
        & =d f_{., \widetilde{\mathcal{H}}}+d f_{., \widetilde{\mathcal{V}}},
    \end{aligned}
\end{equation}
where
\begin{equation}
    \begin{aligned}
        d f_{\mathcal{H}, \widetilde{\mathcal{H}}} & =\pi_{\widetilde{\mathcal{H}}} \circ d f \circ \pi_{\mathcal{H}}, \quad d f_{\mathcal{H}, \widetilde{\mathcal{V}}}=\pi_{\widetilde{\mathcal{V}}} \circ d f \circ \pi_{\mathcal{H}}, \\
        d f_{\mathcal{V}, \widetilde{\mathcal{H}}} & =\pi_{\widetilde{\mathcal{H}}} \circ d f \circ \pi_{\mathcal{V}}, \quad d f_{\mathcal{V}, \widetilde{\mathcal{V}}}=\pi_{\widetilde{\mathcal{V}}} \circ d f \circ \pi_{\mathcal{V}} ,\\
        d f_{\mathcal{H}} & =d f \circ \pi_{\mathcal{H}}, \quad d f_{\mathcal{V}}=d f \circ \pi_{\mathcal{V}}, \\
        d f_{., \widetilde{\mathcal{H}}} & =\pi_{\widetilde{\mathcal{H}}} \circ d f, \quad d f_{., \widetilde{\mathcal{V}}}=\pi_{\widetilde{\mathcal{V}}} \circ d f .
    \end{aligned}
\end{equation}
All these partial differentials are sections of $T^* M \otimes f^{-1} T N$. The bundle $T^* M \otimes$ $f^{-1} T N$ has the induced connection $\nabla \otimes f^{-1} \widetilde{\nabla}$, where $f^{-1} \widetilde{\nabla}$ is the pull-back connection from the Levi-Civita connection $\widetilde{\nabla}$ in $N$. For simplicity, we sometimes write $f^{-1} \widetilde{\nabla}$ as $\widetilde{\nabla}$ when the meaning is clear. Then the horizontal second fundamental form of $f$ with respect to $(\nabla, \widetilde{\nabla})$ is defined by
\begin{equation}
\beta_{\mathcal{H}, \widetilde{\mathcal{H}}}(X, Y)=\widetilde{\nabla}_X d f_{\mathcal{H}, \widetilde{\mathcal{H}}}(Y)-d f_{\mathcal{H}, \widetilde{\mathcal{H}}}\left(\nabla_X Y\right)
\end{equation}
for any $X, Y \in \Gamma(T M)$. There are two special cases : (i) $(M, \mathcal{F})=M$ is a point foliation; (ii) $(N, \widetilde{\mathcal{F}})=N$ is a point foliation. For these two cases, their horizontal second fundamental forms are given by $\beta_{\widetilde{\mathcal{H}}}:=\beta_{TM, \widetilde{\mathcal{H}}}(\cdot, \cdot)$ and $\beta_{\mathcal{H}}:=\beta_{\mathcal{H}, TN}(\cdot, \cdot)$ respectively.

For any $p \in M$, we let $\left\{e_A\right\}_{A=1}^{m+k}$ be an adapted frame field around $p$ and let $\left\{\widetilde{e}_{\widetilde{A}}\right\}_{\widetilde{A}=1}^{n+l}$ be an adapted frame field around $f(p) \in N$. From now on, we shall make use of the following convention on the ranges of indices in $M$ and $N$ respectively:
$$
\begin{aligned}
& 1 \leq A, B, C, \ldots \leq m+k; \\
& 1 \leq i, j, k, \ldots \leq m, \quad m+1 \leq \alpha, \beta, \gamma, \ldots \leq m+k; \\
& 1 \leq \widetilde{A}, \widetilde{B}, \widetilde{C}, \ldots \leq n+l; \\
& 1 \leq \widetilde{i}, \widetilde{j}, \widetilde{k}, \ldots \leq n, \quad n+1 \leq \widetilde{\alpha}, \widetilde{\beta}, \widetilde{\gamma}, \ldots \leq n+l,
\end{aligned}
$$
and we shall agree that repeated indices are summed over the respective ranges. Using the frame fields, we write
\begin{equation}\label{eqdf}
    \begin{aligned}
        d f\left(e_A\right) & =f_A^{\widetilde{B}} \widetilde{e}_{\widetilde{B}}=f_A^{\widetilde{i}} \widetilde{e}_i+f_A^{\widetilde{\alpha}} \widetilde{e}_{\widetilde{\alpha}} ,\\
        d f_{\mathcal{H}, \widetilde{\mathcal{H}}}\left(e_i\right) & =f_i^{\widetilde{j}} \widetilde{e}_{\widetilde{j}}, \quad d f_{\mathcal{H}, \widetilde{\mathcal{V}}}\left(e_i\right)=f_i^{\widetilde{\alpha}} \widetilde{e}_{\widetilde{\alpha}}, \\
        d f_{\mathcal{V}, \widetilde{\mathcal{H}}}\left(e_\alpha\right) & =f_\alpha^{\widetilde{i}} \widetilde{e}_{\widetilde{i}}, \quad d f_{\mathcal{V}, \widetilde{\mathcal{V}}}\left(e_\alpha\right)=f_\alpha^{\widetilde{\beta}} \widetilde{e}_{\widetilde{\beta}} .
    \end{aligned}
\end{equation}
A smooth map $f: M \rightarrow N$ is said to be foliated if $d f(\mathcal{V}) \subset f^{-1} \widetilde{\mathcal{V}}$, or equivalently, $f_\alpha^{\widetilde{i}}=0$ for any $\alpha$ and $\widetilde{i}$. Using $\eqref{eqdf}$, we see that a smooth map $f$ preserves the horizontal distributions, that is, $d f(\mathcal{H}) \subset f^{-1} \widetilde{\mathcal{H}}$ if and only if $f_i^{\widetilde{\alpha}}=0$ for any $i$ and $\widetilde{\alpha}$.

Now we consider the following horizontal energy functional for maps between the foliations:
\begin{equation}\label{eq:hh}
    E_{\mathcal{H}, \widetilde{\mathcal{H}}}(f)=\frac{1}{2} \int_M|d f_{\mathcal{H}, \widetilde{\mathcal{H}}}|^2 d V_g,
\end{equation}
where $d V_g$ is the volume element of the metric $g$, and
\begin{equation}
    |d f_{\mathcal{H}, \widetilde{\mathcal{H}}}|^2=\sum_{i, \widetilde{j}}\left(f_i^{\widetilde{j}}\right)^2.
\end{equation}
We need the following lemma to derive the first variation formula of $E_{\mathcal{H}, \widetilde{\mathcal{H}}}(f)$.

\begin{lemma}[cf. \cite{eellsSelectedTopicsHarmonic1983a}]\label{el}
    Let $f: M \rightarrow N$ be a smooth map. Then
    $$
    \widetilde{\nabla}_X d f(Y)-\widetilde{\nabla}_Y d f(X)-d f([X, Y])=0
    $$
    for any $X, Y \in \Gamma(T M)$.
\end{lemma}

\begin{proposition}\label{1stv}
    Let $f:\left(M^{m+k}, \mathcal{F}^k, g\right) \rightarrow(N^{n+l}, \widetilde{\mathcal{F}}^l, \widetilde{g})$ be a smooth map between two Riemannian foliations and let $\left\{f_t\right\}_{|t|<\varepsilon}$ be a family of maps between $M$ and $N$ with $f_0=f$ and $\left.\frac{\partial f_t}{\partial t}\right|_{t=0}= V \in \Gamma\left(f^{-1} T N\right)$. Then we have
    \begin{equation}
        \frac{d}{d t} E_{\mathcal{H}, \widetilde{\mathcal{H}}}\left(f_t\right)=-\int_M\left\langle V, \tau_{\mathcal{H}, \widetilde{\mathcal{H}}}(f)\right\rangle,
    \end{equation}
    where
    \begin{equation}
                    \begin{aligned}
            \tau_{\mathcal{H}, \widetilde{\mathcal{H}}}(f) & =\operatorname{tr}_g \beta_{\mathcal{H}, \widetilde{\mathcal{H}}}+\operatorname{tr}_g\left(f^* \widetilde{\mathcal{T}}^*+f^* \widetilde{\mathcal{W}}\right) \\
            & =\beta_{\mathcal{H}, \widetilde{\mathcal{H}}}\left(e_i, e_i\right)-d f_{\mathcal{H}, \widetilde{\mathcal{H}}}(\kappa)+\operatorname{tr}_g\left(f^* \widetilde{\mathcal{T}}^*+f^* \widetilde{\mathcal{W}}\right)
        \end{aligned}
    \end{equation}
    and
    \begin{equation}
        \operatorname{tr}_g\left(f^* \widetilde{\mathcal{T}}^*\right)=\widetilde{\mathcal{T}}^*_{d f_{\mathcal{H}, \widetilde{\mathcal{H}}}\left(e_i\right)} d f_{\mathcal{H}, \widetilde{\mathcal{V}}}\left(e_i\right), \quad \operatorname{tr}_g\left(f^* \widetilde{\mathcal{W}}\right)=\widetilde{\mathcal{W}}_{d f_{\mathcal{H}, \widetilde{\mathcal{H}}}\left(e_i\right)} d f_{\mathcal{H}, \widetilde{\mathcal{V}}}\left(e_i\right).
    \end{equation}
\end{proposition}
\begin{proof}
    Let $\Phi: M \times(-\varepsilon, \varepsilon) \rightarrow N$ be the map defined by $\Phi(x, t)=f_t(x)$ for any $(x, t) \in M \times(-\varepsilon, \varepsilon)$. A vector $X \in T M$ may be identified with vector $(X, 0) \in T(M \times(-\varepsilon, \varepsilon))$. This identification gives a corresponding horizontal distribution on $M \times(-\varepsilon, \varepsilon)$, still denoted by $\mathcal{H}$, which is defined by
    $$
    \mathcal{H}_{(x, t)}=\operatorname{span}\left\{(X, 0): X \in \mathcal{H}_x\right\}, \quad(x, t) \in M \times(-\varepsilon, \varepsilon) .
    $$
    
    Let $\left\{e_A\right\}_{A=1}^{m+k}$ be an adapted frame field in $M$. For simplicity, we shall abbreviate the corresponding vector field $\left(e_A, 0\right)$ in $M \times(-\varepsilon, \varepsilon)$ as $e_A$ for each $1 \leq A \leq m+k$ in the following. Applying Lemma \ref{el}, we have
    \begin{equation}\label{eq26}
        \begin{aligned}
            & \frac{d}{d t} E_{\mathcal{H}, \widetilde{\mathcal{H}}}\left(f_t\right) \\
            = & \int_M\left\langle\widetilde{\nabla}_{\frac{\partial}{\partial t}} d \Phi_{\mathcal{H}, \widetilde{\mathcal{H}}}\left(e_i\right), d \Phi_{\mathcal{H}, \widetilde{\mathcal{H}}}\left(e_i\right)\right\rangle d V_g \\
            = & \int_M\left\langle\widetilde{\nabla}_{\frac{\partial}{\partial t}} d \Phi\left(e_i\right)-\widetilde{\nabla}_{\frac{\partial}{\partial t}} d \Phi_{\mathcal{H}, \widetilde{\mathcal{V}}}\left(e_i\right), d \Phi_{\mathcal{H}, \widetilde{\mathcal{H}}}\left(e_i\right)\right\rangle d V_g \\
            = & \int_M\left\{\left\langle\widetilde{\nabla}_{e_i} d \Phi\left(\frac{\partial}{\partial t}\right), d \Phi_{\mathcal{H}, \widetilde{\mathcal{H}}}\left(e_i\right)\right\rangle-\left\langle\widetilde{\nabla}_{\frac{\partial}{\partial t}} d \Phi_{\mathcal{H}, \widetilde{\mathcal{V}}}\left(e_i\right), d \Phi_{\mathcal{H}, \widetilde{\mathcal{H}}}\left(e_i\right)\right\rangle\right\} d V_g.
        \end{aligned}
    \end{equation}
    According to $\eqref{eqdf}$, we may write
    \begin{equation}\label{eq27}
        \begin{aligned}
            d \Phi_{\mathcal{H}, \widetilde{\mathcal{H}}}\left(e_i\right) & =\Phi_i^{\widetilde{j}} \widetilde{e}_{\widetilde{j}}, \quad d \Phi_{\mathcal{H}, \widetilde{\mathcal{V}}}\left(e_i\right)=\Phi_i^{\widetilde{\beta}} \widetilde{e}_{\widetilde{\beta}} \\
            d \Phi\left(\frac{\partial}{\partial t}\right) & =\Phi_0^{\widetilde{i}} \widetilde{e}_{\widetilde{i}}+\Phi_0^{\widetilde{\alpha}} \widetilde{e}_{\widetilde{\alpha}} .
        \end{aligned}
    \end{equation}
    Using $\eqref{eq10}$, $\eqref{eq11}$, $\eqref{eq14}$ and $\eqref{eq27}$, we rewrite the second integrand in $\eqref{eq26}$ as follows:
    \begin{align}
        &\quad \left\langle\widetilde{\nabla}_{\frac{\partial}{\partial t}} d \Phi_{\mathcal{H}, \widetilde{\mathcal{V}}}\left(e_i\right), d \Phi_{\mathcal{H}, \widetilde{\mathcal{H}}}\left(e_i\right)\right\rangle=\left\langle\widetilde{\nabla}_{\frac{\partial}{\partial t}}\left(\Phi_i^{\widetilde{\beta}} \widetilde{e}_{\widetilde{\beta}}\right), d \Phi_{\mathcal{H}, \widetilde{\mathcal{H}}}\left(e_i\right)\right\rangle \notag\\
        & =\Phi_i^{\widetilde{\beta}}\left\langle\widetilde{\nabla}_{\frac{\partial}{\partial t}} \widetilde{e}_{\widetilde{\beta}}, d \Phi_{\mathcal{H}, \widetilde{\mathcal{H}}}\left(e_i\right)\right\rangle=\Phi_i^{\widetilde{\beta}}\left\langle\widetilde{\nabla}_{d \Phi\left(\frac{\partial}{\partial t}\right)} \widetilde{e}_{\widetilde{\beta}}, d \Phi_{\mathcal{H}, \widetilde{\mathcal{H}}}\left(e_i\right)\right\rangle \notag\\
        & =\Phi_i^{\widetilde{\beta}}\left\langle\widetilde{\nabla}_{\Phi_0^{\widetilde{i}}\widetilde{e}_{\widetilde{i}}+\Phi_0^{\widetilde{\alpha}} \widetilde{e}_{\widetilde{\alpha}}}  \widetilde{e}_{\widetilde{\beta}}, d \Phi_{\mathcal{H}, \widetilde{\mathcal{H}}}\left(e_i\right)\right\rangle \notag\\
        & =\Phi_i^{\widetilde{\beta}} \Phi_0^{\widetilde{i}} \Phi_i^{\widetilde{j}}\left\langle\widetilde{\nabla}_{\widetilde{e}_{\widetilde{i}}} \widetilde{e}_{\widetilde{\beta}}, \widetilde{e}_{\widetilde{j}}\right\rangle+\Phi_i^{\widetilde{\beta}} \Phi_0^{\widetilde{\alpha}}\left\langle\widetilde{\nabla}_{\widetilde{e}_{\widetilde{\alpha}}} \widetilde{e}_{\widetilde{\beta}}, d \Phi_{\mathcal{H}, \widetilde{\mathcal{H}}}\left(e_i\right)\right\rangle \notag\\
        & =-\Phi_i^{\widetilde{\beta}} \Phi_0^{\widetilde{i}} \Phi_i^{\widetilde{j}}\left\langle\widetilde{e}_{\widetilde{\beta}}, \widetilde{\nabla}_{\widetilde{e}_i} \widetilde{e}_{\widetilde{j}}\right\rangle+\left\langle\mathcal{S}\left(\Phi_0^{\widetilde{\alpha}} \widetilde{e}_{\widetilde{\alpha}}, \Phi_i^{\widetilde{\beta}} \widetilde{e}_{\widetilde{\beta}}\right), d \Phi_{\mathcal{H}, \widetilde{\mathcal{H}}}\left(e_i\right)\right\rangle \notag\\
        & =\left\langle\Phi_i^{\widetilde{\beta}} \widetilde{e}_{\widetilde{\beta}}, \widetilde{\mathcal{T}}_{\Phi_i^{\widetilde{j}} \widetilde{e}_{\widetilde{j}}} \Phi_0^{\widetilde{i}} \widetilde{e}_{\widetilde{i}}\right\rangle+\left\langle d \Phi_{\cdot, \widetilde{\mathcal{V}}}\left(\frac{\partial}{\partial t}\right), \widetilde{\mathcal{W}}_{d \Phi_{\mathcal{H}, \widetilde{\mathcal{H}}}\left(e_i\right)} d \Phi_{\mathcal{H}, \widetilde{\mathcal{V}}}\left(e_i\right)\right\rangle \notag\\
        & =\left\langle d \Phi_{\cdot, \widetilde{\mathcal{H}}}\left(\frac{\partial}{\partial t}\right), \widetilde{\mathcal{T}}_{d \Phi_{\mathcal{H}, \widetilde{\mathcal{H}}}\left(e_i\right)}^* d \Phi_{\mathcal{H}, \widetilde{\mathcal{V}}}\left(e_i\right)\right\rangle+\left\langle d \Phi_{\cdot, \widetilde{\mathcal{V}}}\left(\frac{\partial}{\partial t}\right), \widetilde{\mathcal{W}}_{d \Phi_{\mathcal{H}, \widetilde{\mathcal{H}}}\left(e_i\right)} d \Phi_{\mathcal{H}, \widetilde{\mathcal{V}}}\left(e_i\right)\right\rangle \notag\\
        & =\left\langle d \Phi\left(\frac{\partial}{\partial t}\right), \widetilde{\mathcal{T}}_{d \Phi_{\mathcal{H}, \widetilde{\mathcal{H}}}\left(e_i\right)}^* d \Phi_{\mathcal{H}, \widetilde{\mathcal{V}}}\left(e_i\right)+\widetilde{\mathcal{W}}_{d \Phi_{\mathcal{H}, \widetilde{\mathcal{H}}}\left(e_i\right)} d \Phi_{\mathcal{H}, \widetilde{\mathcal{V}}}\left(e_i\right)\right\rangle . \label{eq28}
    \end{align}
    Set
    \begin{equation}\label{eq29}
        \operatorname{tr}_g(f^* \widetilde{\mathcal{T}}^*)=\widetilde{\mathcal{T}}^{*}_{d f_{\mathcal{H}, \widetilde{\mathcal{H}}}\left(e_i\right)} d f_{\mathcal{H}, \widetilde{\mathcal{V}}}\left(e_i\right), \operatorname{tr}_g(f^* \widetilde{\mathcal{W}})=\widetilde{\mathcal{W}}_{d f_{\mathcal{H}, \widetilde{\mathcal{H}}}\left(e_i\right)} d f_{\mathcal{H}, \widetilde{\mathcal{V}}}\left(e_i\right).
    \end{equation}
    From $\eqref{eq26}$, $\eqref{eq28}$ and $\eqref{eq29}$, we obtain
    \begin{equation}\label{eq:320}
        \begin{aligned}
            & \left.\frac{d}{d t} E_{\mathcal{H}, \widetilde{\mathcal{H}}}\left(f_t\right)\right|_{t=0} \\
            = & \int_M\left\{\left\langle\widetilde{\nabla}_{e_i}  V, d \Phi_{\mathcal{H}, \widetilde{\mathcal{H}}}\left(e_i\right)\right\rangle-\left\langle V, \operatorname{tr}_g\left(f^* \widetilde{\mathcal{T}}^*\right)+\operatorname{tr}_g\left(f^* \widetilde{\mathcal{W}}\right)\right\rangle\right\} d V_g \\
            = & \int_M\left\{e_i\left\langle V, d f_{\mathcal{H}, \widetilde{\mathcal{H}}}\left(e_i\right)\right\rangle-\left\langle V, \widetilde{\nabla}_{e_i} d f_{\mathcal{H}, \widetilde{\mathcal{H}}}\left(e_i\right)\right\rangle-\left\langle V, \operatorname{tr}_g\left(f^* \widetilde{\mathcal{T}}^*\right)+\operatorname{tr}_g\left(f^* \widetilde{\mathcal{W}}\right)\right\rangle\right\} d V_g .
        \end{aligned}
    \end{equation}
    Define a 1-form $\theta$ on $M$ by
    $$
    \theta(X)=\left\langle V, d f_{\mathcal{H}, \widetilde{\mathcal{H}}}(X)\right\rangle, \quad \forall X \in T M .
    $$
    The codifferential of $\theta$ is given by
    \begin{equation}\label{eq:330}
        \begin{aligned}
            \delta \theta & =-\left(\nabla_{e_A} \theta\right)\left(e_A\right) \\
            & =-e_A \theta\left(e_A\right)+\theta\left(\nabla_{e_A} e_A\right) \\
            & =-e_i\left\langle V, d f_{\mathcal{H}, \widetilde{\mathcal{H}}}\left(e_i\right)\right\rangle+\left\langle V, d f_{\mathcal{H}, \widetilde{\mathcal{H}}}\left(\nabla_{e_A} e_A\right)\right\rangle \\
            & =-e_i\left\langle V, d f_{\mathcal{H}, \widetilde{\mathcal{H}}}\left(e_i\right)\right\rangle+\left\langle V, d f_{\mathcal{H}, \widetilde{\mathcal{H}}}\left(\nabla_{e_i} e_i\right)\right\rangle+\left\langle V, d f_{\mathcal{H}, \widetilde{\mathcal{H}}}(\kappa)\right\rangle .
        \end{aligned}
    \end{equation}
    From $\eqref{eq:320}$, $\eqref{eq:330}$ and Stokes' theorem, we conclude that
    \begin{equation*}
        \left.\frac{d}{d t} E_{\mathcal{H}, \widetilde{\mathcal{H}}}\left(f_t\right)\right|_{t=0}=-\int_M\left\langle V, \operatorname{tr}_g \beta_{\mathcal{H}, \widetilde{\mathcal{H}}}+\operatorname{tr}_g(f^* \widetilde{\mathcal{T}}^*)+\operatorname{tr}_g(f^* \widetilde{\mathcal{W}})\right\rangle d V_g
    \end{equation*}
    and
    \begin{equation}
        \operatorname{tr}_g \beta_{\mathcal{H}, \widetilde{\mathcal{H}}}=\beta_{\mathcal{H}, \widetilde{\mathcal{H}}}\left(e_i, e_i\right)-d f(\kappa) .
    \end{equation}
    This proves Proposition \ref{1stv}.
\end{proof}

\begin{remark}\label{re:tension}
Let $\{e_{A}\}$ and $\{\widetilde{e}_{\widetilde{A}}\}$ be adapted
frame fields around $p\in M$ and $q=f(p)\in N$ respectively, such that $\pi
_{\mathcal{H}}(\nabla _{e_{i}}e_{j})_{p}=0$, $\pi _{\mathcal{H}}([e_{\alpha },e_{i}])=0
$, $\pi _{\widetilde{\mathcal{H}}}(\widetilde{\nabla }_{\widetilde{e}_{
\widetilde{i}}}\widetilde{e}_{\widetilde{j}})_{q}=0$ and $\pi _{\widetilde{\mathcal{H}}%
}([\widetilde{e}_{\widetilde{\alpha }},\widetilde{e}_{\widetilde{i}}])=0$
(see \S \ref{se:Riemfoli}). Using $\eqref{wten}, \eqref{ttensor}$ and Lemma \ref{le:td}, we derive
\begin{eqnarray}\label{eq:35}
\tau _{\mathcal{H},\widetilde{\mathcal{H}}}(f)_{p} &=&\beta _{\mathcal{H},%
\widetilde{\mathcal{H}}}(e_{i},e_{i})-df_{\mathcal{H},\widetilde{\mathcal{H}}%
}(\kappa )+\widetilde{\mathcal{T}}_{df_{\mathcal{H},\widetilde{\mathcal{H}}%
}(e_{i})}^{\ast }df_{\mathcal{H},\widetilde{\mathcal{V}}}(e_{i})+\widetilde{%
\mathcal{W}}_{df_{\mathcal{H},\widetilde{\mathcal{H}}}(e_{i})}df_{\mathcal{H}%
,\widetilde{\mathcal{V}}}(e_{i})  \nonumber \\
&=&\widetilde{\nabla }_{e_{i}}(f_{i}^{\widetilde{j}}\widetilde{e}_{%
\widetilde{j}})-df_{\mathcal{H},\widetilde{\mathcal{H}}}(\kappa )+f_{i}^{%
\widetilde{j}}f_{i}^{\widetilde{\beta }}\widetilde{\mathcal{T}}_{\widetilde{e%
}_{\widetilde{j}}}^{\ast }\widetilde{e}_{\widetilde{\beta }}+f_{i}^{%
\widetilde{j}}f_{i}^{\widetilde{\beta }}\widetilde{\mathcal{W}}_{\widetilde{e%
}_{\widetilde{j}}}\widetilde{e}_{\widetilde{\beta }}  \nonumber \\
&=&e_{i}(f_{i}^{\widetilde{j}})\widetilde{e}_{\widetilde{j}}+f_{i}^{%
\widetilde{j}}f_{i}^{\widetilde{A}}\widetilde{\nabla }_{\widetilde{e}_{%
\widetilde{A}}}\widetilde{e}_{\widetilde{j}}-df_{\mathcal{H},\widetilde{%
\mathcal{H}}}(\kappa )-f_{i}^{\widetilde{j}}f_{i}^{\widetilde{\beta }}%
\widetilde{\nabla }_{\widetilde{e}_{\widetilde{\beta }}}\widetilde{e}_{%
\widetilde{j}}  \nonumber \\
&=&e_{i}(f_{i}^{\widetilde{j}})\widetilde{e}_{\widetilde{j}}-df_{\mathcal{H},%
\widetilde{\mathcal{H}}}(\kappa )+f_{i}^{\widetilde{j}}f_{i}^{\widetilde{k}%
}\pi _{\widetilde{\mathcal{V}}}(\widetilde{\nabla }_{\widetilde{e}_{%
\widetilde{k}}}\widetilde{e}_{\widetilde{j}})  \nonumber \\
&=&e_{i}(f_{i}^{\widetilde{j}})\widetilde{e}_{\widetilde{j}}-df_{\mathcal{H},%
\widetilde{\mathcal{H}}}(\kappa ).  
\end{eqnarray}%
Therefore this proves that $\tau _{\mathcal{H},\widetilde{\mathcal{H}}%
}(f)\in \Gamma (f^{-1}\widetilde{\mathcal{H}})$. Notice that if $f$ is foliated, then
the right hand side of $\eqref{eq:35}$ is exactly the transverse tension field $\tau
_{T}(f)$ of the energy functional $E_{\mathcal{H},\widetilde{\mathcal{H}}}(f)
$ under foliated variations considered in \cite{barlettaTransversallyHolomorphicMaps1998a}, \cite{chenBochnertypeFormulasTransversally2012a} and
\cite{dragomirHarmonicMapsFoliated2013a}.
Furthermore, if $\mathcal{F}$ is harmonic, then $\tau _{\mathcal{H},%
\widetilde{\mathcal{H}}}(f)=0$ if and only if $\tau (\bar{f})=0$, where $\bar{f}$ is the locally induced map of $f$ and $\tau(\bar{f})$ is the usual tension field of $\bar{f}$.
\end{remark}

\begin{definition}
    A smooth map $f:(M^{m+k}, \mathcal{F}^k, g) \rightarrow(N^{n+l}, \widetilde{\mathcal{F}}^l, \widetilde{g})$ between two Riemannian foliations is referred to as a horizontally harmonic map if it is a critical map of $E_{\mathcal{H}, \widetilde{\mathcal{H}}}(f)$, that is,
    \begin{equation}
        \operatorname{tr}_g \beta_{\mathcal{H}, \widetilde{\mathcal{H}}}+\operatorname{tr}_g\left(f^* \widetilde{\mathcal{T}}^*+f^* \widetilde{\mathcal{W}}\right)=0.
    \end{equation}
    In particular, if $f$ is foliated, then $f$ is called a transversally harmonic map (or equivalently, a $(\mathcal{F},\widetilde{\mathcal{F}})$-harmonic map
    according to \cite{dragomirHarmonicMapsFoliated2013a}).
\end{definition}

\begin{remark} 
    As we have mentioned in  the Introduction, foliated critical
    points of $E_{\mathcal{H},\widetilde{\mathcal{H}}}(f)$, called $(\mathcal{F},\widetilde{\mathcal{F}})$-harmonic maps in \cite{dragomirHarmonicMapsFoliated2013a}, are different from the
    transversally harmonic maps defined in \cite{konderakTransversallyHarmonicMaps2003}, and the notion of
    horizontally harmonic maps is a slight generalization of the notions 
    defined in \cite{chenBochnertypeFormulasTransversally2012a} and  \cite{dragomirHarmonicMapsFoliated2013a} (see also Remark \ref{re:tension}).
    Notice that a horizontal energy functional similar to $\eqref{eq:hh}$ also appeared in
    \cite{petitMokSiuYeungTypeFormulas2009} for maps between pseudo-Hermitian manifolds that are not necessarily
    Riemannian foliations (see also \cite{dong$HwidetildeH$harmonicMapsPseudoHermitian2016a}). 
\end{remark}

\begin{corollary}
    Let $f:\left(M^{m+k}, \mathcal{F}^k, g\right) \rightarrow\left(N^{n+l}, \widetilde{g}\right)$ be a map from a Riemannian foliation to a Riemannian manifold (that is, a Riemannian point foliation). Then $f$ is a horizontally harmonic map if and only if
    $$
    \operatorname{tr}_g \beta_{\mathcal{H}}=\beta_{\mathcal{H}}\left(e_i, e_i\right)-d f(\kappa)=0 .
    $$
\end{corollary}

\begin{remark}
     When $\mathcal{H}$ satisfies the Hörmander's condition, a horizontally harmonic map $f:\left(M^{m+k}, \mathcal{F}^k, g\right) \rightarrow\left(N^{n+l}, \widetilde{g}\right)$ is also called a subelliptic harmonic map (cf. \cite{dongEellsSampsonType2021a}). 
     We can expect subelliptic harmonic maps to have better analytic properties than general horizontally harmonic maps.
\end{remark}

\begin{proposition}
    Let $f:(M,g) \to (N,h)$ be a smooth map between Riemannian manifolds. Then $f$ is harmonic if and only if the differential map $f_{*}:(TM,\widehat{g}) \to (TN, \widehat{h})$ is horizontally harmonic.
\end{proposition}
\begin{proof}
    Note that $f_{*}$ is foliated and the leaves of $(TM,\widehat{g})$ are totally geodesic (and therefore minimal). Applying Remark \ref{re:tension}, we know that $f_{*}$ is horizontally harmonic if and only if $f$ is harmonic.
\end{proof}
\begin{remark}
    From the above proposition, we can get many examples of horizontally harmonic maps between the tangent bundles from harmonic maps between the base Riemannian manifolds. 
\end{remark}

A trivial horizontally harmonic map is the map $f: M \rightarrow N$ with $d f_{\mathcal{H}, \widetilde{\mathcal{H}}} \equiv 0$. Such a map is said to be horizontally constant. Let us check the geometric meaning of this property.

\begin{lemma}
    Let $f:(M^{m+k}, \mathcal{F}^k, g) \rightarrow(N^{n+l}, \widetilde{\mathcal{F}}^l, \widetilde{g})$ be a map with $d f_{\mathcal{H}, \widetilde{\mathcal{H}}} \equiv 0$. Then for any horizontal curve $\gamma, f(\gamma)$ is contained in a single leaf of $\widetilde{\mathcal{F}}^l$. In particular, (i) if $\mathcal{H}$ satisfies Hörmander's condition, then $f(M)$ is contained in a single leaf; (ii) if $\mathcal{H}$ satisfies $\mathcal{T} \equiv 0$, then for each connected integral submanifold $P$ of $\mathcal{H}, f(P)$ is contained in a single leaf.
\end{lemma}
\begin{proof}
     Let $p \in M$ and $q=f(p) \in N$. Choose foliated coordinate charts $U$ and $\widetilde{U}$ around $p$ and $q$, respectively, such that $f(U) \subset \widetilde{U}$. The coordinate chart $\widetilde{U}$ induces a natural projection $\widetilde{\pi}: \widetilde{U} \rightarrow \widetilde{\mathcal{B}}_{\widetilde{U}}$ to a local base manifold. Set $\widehat{f}=\widetilde{\pi} \circ f$. The condition $d f_{\mathcal{H}, \widetilde{\mathcal{H}}} \equiv 0$ implies that $d \widehat{f}_{\mathcal{H}} \equiv 0$, that is, $\widehat{f}$ is constant along any horizontal curve in $U$. This means $f$ maps any horizontal curve in $U$ into a leaf in $\widetilde{U}$. Let $\gamma:[0, l] \rightarrow M$ be any horizontal curve in $M$. We may cover $\gamma$ by a finite number of coordinate charts. Then it is easy to see that $f(\gamma)$ is contained in a single leaf too.
    
    First, we assume that $\mathcal{H}$ satisfies the Hörmander's condition. Since any two points can be joined by a horizontal curve, we find that $f(M)$ is contained in a single leaf. This proves (i).
    
    Next, we suppose that $\mathcal{H}$ satisfies $\mathcal{T} \equiv 0$, that is, $\mathcal{H}$ is integrable. Let $P$ be any connected integral submanifold of $\mathcal{H}$ and $p \in P$. Then for any $q \in P$, there is a horizontal curve $\gamma$ joining $p$ and $q$. By the previous argument, we see that $f(\gamma)$ is contained in a single leaf. Since $p$ is fixed and $q$ is arbitrary, we conclude that $f(P)$ is contained in a single leaf. Then (ii) is proved.
\end{proof}

\begin{remark}
    In case (ii) of the above lemma, let $\widehat{\mathcal{F}}$ denote the foliation consisting of the integral submanifolds of $\mathcal{H}$. Then the result means that $f$ : $(M^{m+k}, \widehat{\mathcal{F}}^m) \rightarrow(N^{n+l}, \widetilde{\mathcal{F}}^l)$ is foliated.
\end{remark}
~

\section{The stress-energy tensor and conservation laws}\label{se:streeenergy}

In this section, we introduce the stress-energy tensor for maps between two Riemannian foliations, and then investigate the conservation laws for the critical maps of the energy functional $E_{\mathcal{H}, \widetilde{\mathcal{H}}}(f)$.

Let us first recall briefly the notion of the stress-energy tensor defined on vector bundle valued $p$-forms (cf. \cite{dongVanishingTheoremsVector2011a} for details). Let $\xi: E \rightarrow M$ be a Riemannian vector bundle over a Riemannian manifold $(M, g)$. Set $\mathcal{A}^p(\xi)=$ $\Gamma\left(\Lambda^p T^* M \otimes \xi\right)$. We consider the following energy functional
\begin{equation}
    \mathcal{E}(\omega)=\frac{1}{2} \int_M|\omega|^2 d V_g
\end{equation}
for any $\omega \in \mathcal{A}^p(\xi)$. Then the stress-energy tensor associated with the $\mathcal{E}$-energy functional is given by:
\begin{equation}
    S_\omega(X, Y)=\frac{|\omega|^2}{2} g(X, Y)-(\omega \odot \omega)(X, Y)
\end{equation}
for any $X, Y \in T M$, where $(\omega \odot \omega)$ denotes a 2-tensor defined by
\begin{equation}\label{2ten}
    (\omega \odot \omega)(X, Y)=\left\langle i_X \omega, i_Y \omega\right\rangle .
\end{equation}
Here $i_X$ is the interior product with respect to $X$. It is known that $S_\omega$ is a useful tool for studying the $\mathcal{E}$-energy functional.

Noting that $E_{\HH}(f)=\mathcal{E}(\omega)$ with $\omega=d f_{\mathcal{H}, \widetilde{\mathcal{H}}}$, we obtain the following stress-energy tensor associated with the energy functional $E_{\mathcal{H}, \widetilde{\mathcal{H}}}(f)$:
\begin{equation}\label{stressen}
    S_{\mathcal{H}, \widetilde{\mathcal{H}}}(f)=\frac{|d f_{\mathcal{H}, \widetilde{\mathcal{H}}}|^2}{2} g-d f_{\mathcal{H}, \widetilde{\mathcal{H}}} \odot d f_{\mathcal{H}, \widetilde{\mathcal{H}}}.
\end{equation}
According to $\eqref{2ten}$, we have
\[(d f_{\mathcal{H}, \widetilde{\mathcal{H}}} \odot d f_{\mathcal{H}, \widetilde{\mathcal{H}}})(X, Y)=\langle d f_{\mathcal{H}, \widetilde{\mathcal{H}}}(X), d f_{\mathcal{H}, \widetilde{\mathcal{H}}}(Y)\rangle,\]
for $X, Y \in T M$. As a 2-tensor field, the divergence of $S_{\mathcal{H}, \widetilde{\mathcal{H}}}(f)$ is a 1-form on $M$, defined by
\begin{equation}
    (\operatorname{div} S_{\mathcal{H}, \widetilde{\mathcal{H}}}(f))(X)=(\nabla_{e_A} S_{\mathcal{H}, \widetilde{\mathcal{H}}}(f))\left(e_A, X\right)
\end{equation}
for any $X \in \mathfrak{X}(M)$, where $\left\{e_A\right\}$ is any local orthonormal frame field of $(M, g)$.

\begin{theorem}
 Let $f:\left(M^{m+k}, \mathcal{F}^k , g\right) \rightarrow(N^{n+l}, \widetilde{\mathcal{F}}^l , \widetilde{g})$ be a smooth map between two Riemannian foliations and let $S_{\mathcal{H}, \widetilde{\mathcal{H}}}(f)$ be the stress-energy tensor defined by $\eqref{stressen}$. Then
    $$
    \begin{aligned}
        & (\operatorname{div} S_{\mathcal{H}, \widetilde{\mathcal{H}}}(f))(X) \\
        = & -\left\langle\tau_{\mathcal{H}, \widetilde{\mathcal{H}}}(f), d f(X)\right\rangle+\left\langle  \operatorname{tr}_g \beta_{\mathcal{H}, \widetilde{\mathcal{H}}}(f), d f_{\mathcal{V}, \widetilde{\mathcal{H}}}(X)\right\rangle+\left\langle d f_{\mathcal{H}, \widetilde{\mathcal{H}}}\left(e_i\right),\left(\widetilde{\nabla}_{e_i} d f_{\mathcal{V}, \widetilde{\mathcal{H}}}\right)(X)\right\rangle \\
        & +\left\langle d f_{\mathcal{V}, \widetilde{\mathcal{H}}}\left(\mathcal{T}_X e_i\right)-d f_{\mathcal{V}, \widetilde{\mathcal{H}}}\left(\mathcal{W}_{e_i}(X)\right), d f_{\mathcal{H}, \widetilde{\mathcal{H}}}\left(e_i\right)\right\rangle
    \end{aligned}
    $$
    for any $X \in T M$. 
\end{theorem}
\begin{proof}
    Let $p \in M$ and $\left\{e_A\right\}_{A=1}^{m+k}$ be any adapted frame field around $p$ such that $e_i(1 \leq i \leq m)$ is basic, and
    \begin{equation}\label{frame}
        \pi_{\mathcal{H}}\left(\nabla_Z e_i\right)_p=0, \quad 1 \leq i \leq m,
    \end{equation}
    for any $Z \in \mathfrak{X}_{\mathcal{H}}$. For any $X \in \mathfrak{X}(M)$, we compute the divergence of $S_{\mathcal{H}, \widetilde{\mathcal{H}}}(f)$ at $p$ as follows:
    \begin{align}
        \left(\operatorname{div} S_{\mathcal{H}, \widetilde{\mathcal{H}}}(f)\right)(X)= & \left(\nabla_{e_A} S_{\mathcal{H}, \widetilde{\mathcal{H}}}(f)\right)\left(e_A, X\right) \notag\\
        = & e_A\left(S_{\mathcal{H}, \widetilde{\mathcal{H}}}(f)\left(e_A, X\right)\right)-S_{\mathcal{H}, \widetilde{\mathcal{H}}}(f)\left(\nabla_{e_A} e_A, X\right) \notag\\
        & -S_{\mathcal{H}, \widetilde{\mathcal{H}}}(f)\left(e_A, \nabla_{e_A} X\right) \notag\\
        = & e_A[\frac{|d f_{\mathcal{H}, \widetilde{\mathcal{H}}}|^2}{2}\left\langle e_A, X\right\rangle]-e_A\left\langle d f_{\mathcal{H}, \widetilde{\mathcal{H}}}\left(e_A\right), d f_{\mathcal{H}, \widetilde{\mathcal{H}}}(X)\right\rangle \notag\\
        & -\frac{\left|d f_{\mathcal{H}, \widetilde{\mathcal{H}}}\right|^2}{2}\left\langle\nabla_{e_A} e_A, X\right\rangle+\left\langle d f_{\mathcal{H}, \widetilde{\mathcal{H}}}\left(\nabla_{e_A} e_A\right), d f_{\mathcal{H}, \widetilde{\mathcal{H}}}(X)\right\rangle \notag\\
        & -\frac{\left|d f_{\mathcal{H}, \widetilde{\mathcal{H}}}\right|^2}{2}\left\langle e_A, \nabla_{e_A} X\right\rangle+\left\langle d f_{\mathcal{H}, \widetilde{\mathcal{H}}}\left(e_A\right), d f_{\mathcal{H}, \widetilde{\mathcal{H}}}\left(\nabla_{e_A} X\right)\right\rangle \notag\\
        = & e_A(\frac{|d f_{\mathcal{H}, \widetilde{\mathcal{H}}}|^2}{2})\left\langle e_A, X\right\rangle-\left\langle\left(\widetilde{\nabla}_{e_A} d f_{\mathcal{H}, \widetilde{\mathcal{H}}}\right)\left(e_A\right), d f_{\mathcal{H}, \widetilde{\mathcal{H}}}(X)\right\rangle \notag\\
        & -\left\langle d f_{\mathcal{H}, \widetilde{\mathcal{H}}}\left(e_A\right), \widetilde{\nabla}_{e_A} d f_{\mathcal{H}, \widetilde{\mathcal{H}}}(X)\right\rangle+\left\langle d f_{\mathcal{H}, \widetilde{\mathcal{H}}}\left(e_A\right), d f_{\mathcal{H}, \widetilde{\mathcal{H}}}\left(\nabla_{e_A} X\right)\right\rangle \notag\\
        =&-\left\langle\left(\widetilde{\nabla}_{e_A} d f_{\mathcal{H}, \widetilde{\mathcal{H}}}\right)\left(e_A\right), d f_{\mathcal{H}, \widetilde{\mathcal{H}}}(X)\right\rangle+\left\langle\widetilde{\nabla}_X d f_{\mathcal{H}, \widetilde{\mathcal{H}}}\left(e_A\right), d f_{\mathcal{H}, \widetilde{\mathcal{H}}}\left(e_A\right)\right\rangle \notag\\
        & -\left\langle d f_{\mathcal{H}, \widetilde{\mathcal{H}}}\left(e_A\right), \widetilde{\nabla}_{e_A} d f_{\mathcal{H}, \widetilde{\mathcal{H}}}(X)\right\rangle+\left\langle d f_{\mathcal{H}, \widetilde{\mathcal{H}}}\left(e_A\right), d f_{\mathcal{H}, \widetilde{\mathcal{H}}}\left(\nabla_{e_A} X\right)\right\rangle \notag\\
         =&-\underbrace{\left\langle\operatorname{tr}_g \beta_{\mathcal{H}, \widetilde{\mathcal{H}}}(f), d f_{\mathcal{H}, \widetilde{\mathcal{H}}}(X)\right\rangle}_{(I)}+\underbrace{\left\langle\widetilde{\nabla}_X d f_{\mathcal{H}, \widetilde{\mathcal{H}}}\left(e_A\right), d f_{\mathcal{H}, \widetilde{\mathcal{H}}}\left(e_A\right)\right\rangle}_{(I I)} \notag\\
        & -\underbrace{\left\langle d f_{\mathcal{H}, \widetilde{\mathcal{H}}}\left(e_A\right), \widetilde{\nabla}_{e_A} d f_{\mathcal{H}, \widetilde{\mathcal{H}}}(X)\right\rangle}_{(I I I)}+\underbrace{\left\langle d f_{\mathcal{H}, \widetilde{\mathcal{H}}}\left(e_A\right), d f_{\mathcal{H}, \widetilde{\mathcal{H}}}\left(\nabla_{e_A} X\right)\right\rangle}_{(I V)} . \label{divs} 
    \end{align}
By $\eqref{dfde}$, we can write the term $(I)$ of $\eqref{divs}$ as 
    \begin{equation}
        \begin{aligned}
            (I) & =\left\langle\operatorname{tr}_g \beta_{\mathcal{H}, \widetilde{\mathcal{H}}}(f), d f(X)-d f_{\mathcal{V}, \widetilde{\mathcal{H}}}(X)-d f_{\mathcal{H}, \widetilde{\mathcal{V}}}(X)-d f_{\mathcal{V}, \widetilde{\mathcal{V}}}(X)\right\rangle \\
            & =\underbrace{\left\langle\operatorname{tr}_g \beta_{\mathcal{H}, \widetilde{\mathcal{H}}}(f), d f(X)\right\rangle}_{(I)\text{-}1}-\underbrace{\left\langle\operatorname{tr}_g \beta_{\mathcal{H}, \widetilde{\mathcal{H}}}(f), d f_{\mathcal{V}, \widetilde{\mathcal{H}}}(X)+d f_{\mathcal{H}, \widetilde{\mathcal{V}}}(X)+d f_{\mathcal{V}, \widetilde{\mathcal{V}}}(X)\right\rangle}_{(I)\text{-}2} .
        \end{aligned}
    \end{equation}
Next, a direct calculation by using $\eqref{dfde}$ and Lemma \ref{el} shows that
        \begin{align}
            (I I)= & \left\langle\widetilde{\nabla}_X d f_{\mathcal{H}, \widetilde{\mathcal{H}}}\left(e_i\right), d f_{\mathcal{H}, \widetilde{\mathcal{H}}}\left(e_i\right)\right\rangle \notag\\
            = & \left\langle\widetilde{\nabla}_X d f\left(e_i\right), d f_{\mathcal{H}, \widetilde{\mathcal{H}}}\left(e_i\right)\right\rangle-\left\langle\widetilde{\nabla}_X d f_{\mathcal{H}, \widetilde{\mathcal{V}}}\left(e_i\right), d f_{\mathcal{H}, \widetilde{\mathcal{H}}}\left(e_i\right)\right\rangle \notag\\
            = & \left\langle\widetilde{\nabla}_{e_i} d f(X), d f_{\mathcal{H}, \widetilde{\mathcal{H}}}\left(e_i\right)\right\rangle+\left\langle d f\left(\left[X, e_i\right]\right), d f_{\mathcal{H}, \widetilde{\mathcal{H}}}\left(e_i\right)\right\rangle \notag\\
            & -\left\langle\widetilde{\nabla}_X d f_{\mathcal{H}, \widetilde{\mathcal{V}}}\left(e_i\right), d f_{\mathcal{H}, \widetilde{\mathcal{H}}}\left(e_i\right)\right\rangle \notag\\
            = & \underbrace{\left\langle\widetilde{\nabla}_{e_i} d f(X), d f_{\mathcal{H}, \widetilde{\mathcal{H}}}\left(e_i\right)\right\rangle}_{(I I)\text{-}1}+\underbrace{\left\langle d f\left(\nabla_X e_i\right), d f_{\mathcal{H}, \widetilde{\mathcal{H}}}\left(e_i\right)\right\rangle}_{(I I)-2} \notag\\
            & -\underbrace{\left\langle d f\left(\nabla_{e_i} X\right), d f_{\mathcal{H}, \widetilde{\mathcal{H}}}\left(e_i\right)\right\rangle}_{(I I)\text{-}3}-\underbrace{\left\langle\widetilde{\nabla}_X d f_{\mathcal{H}, \widetilde{\mathcal{V}}}\left(e_i\right), d f_{\mathcal{H}, \widetilde{\mathcal{H}}}\left(e_i\right)\right\rangle}_{(I I)\text{-}4} .
        \end{align}
In view of 	$\eqref{wten}$, $\eqref{ttensor}$, $\eqref{frame}$ and Lemma \ref{le:td}, the term $(II)\text{-}2$ can be converted into
    \begin{align}
        (I I)\text{-}2 & =\left\langle d f_{\mathcal{H}, \widetilde{\mathcal{H}}}\left(\nabla_X e_i\right)+d f_{\mathcal{V}, \widetilde{\mathcal{H}}}\left(\nabla_X e_i\right), d f_{\mathcal{H}, \widetilde{\mathcal{H}}}\left(e_i\right)\right\rangle \notag\\
        & =\left\langle d f_{\mathcal{H}, \widetilde{\mathcal{H}}}\left(\nabla_{\pi_{\mathcal{V}}(X)} e_i\right)+d f_{\mathcal{V}, \widetilde{\mathcal{H}}}\left(\nabla_{\pi_{\mathcal{H}}(X)} e_i\right)+d f_{\mathcal{V}, \widetilde{\mathcal{H}}}\left(\nabla_{\pi_{\mathcal{V}}(X)} e_i\right), d f_{\mathcal{H}, \widetilde{\mathcal{H}}}\left(e_i\right)\right\rangle \notag\\
        & =\left\langle-d f_{\mathcal{H}, \widetilde{\mathcal{H}}}\left(\mathcal{T}_{e_i}^*\left(\pi_{\mathcal{V}}(X)\right)\right)+d f_{\mathcal{V}, \widetilde{\mathcal{H}}}\left(\mathcal{T}_{\pi_{\mathcal{H}}(X)} e_i\right)-d f_{\mathcal{V}, \widetilde{\mathcal{H}}}\left(\mathcal{W}_{e_i}\left(\pi_{\mathcal{V}}(X)\right)\right), d f_{\mathcal{H}, \widetilde{\mathcal{H}}}\left(e_i\right)\right\rangle \notag\\
        & =\left\langle-d f_{\mathcal{H}, \widetilde{\mathcal{H}}}\left(\mathcal{T}_{e_i}^*(X)\right)+d f_{\mathcal{V}, \widetilde{\mathcal{H}}}\left(\mathcal{T}_X e_i\right)-d f_{\mathcal{V}, \widetilde{\mathcal{H}}}\left(\mathcal{W}_{e_i}(X)\right), d f_{\mathcal{H}, \widetilde{\mathcal{H}}}\left(e_i\right)\right\rangle .
    \end{align}
For the term $(II)\text{-}4$, we deduce in terms of $\eqref{ttensor}$, $\eqref{eq10}$ and $\eqref{eq11}$ that	
\begin{align}
    (I I)\text{-}4 & =\left\langle\widetilde{\nabla}_{{d f_{\cdot, \widetilde{\mathcal{H}}}}(X)+d f_{\cdot, \widetilde{\mathcal{V}}}(X)} d f_{\mathcal{H}, \widetilde{\mathcal{V}}}\left(e_i\right), d f_{\mathcal{H}, \widetilde{\mathcal{H}}}\left(e_i\right)\right\rangle \notag\\
    & =-\left\langle d f_{\mathcal{H}, \widetilde{\mathcal{V}}}\left(e_i\right), \widetilde{\nabla}_{d f_{, \widetilde{\mathcal{H}}}(X)} d f_{\mathcal{H}, \widetilde{\mathcal{H}}}\left(e_i\right)\right\rangle+\left\langle\widetilde{S}\left(d f_{\cdot, \widetilde{\mathcal{V}}}(X), d f_{\mathcal{H}, \widetilde{\mathcal{V}}}\left(e_i\right)\right), d f_{\mathcal{H}, \widetilde{\mathcal{H}}}\left(e_i\right)\right\rangle \notag\\
    & =-\left\langle d f_{\mathcal{H}, \widetilde{\mathcal{V}}}\left(e_i\right), \widetilde{\mathcal{T}}_{d f_{, \widetilde{\mathcal{H}}}(X)} d f_{\mathcal{H}, \widetilde{\mathcal{H}}}\left(e_i\right)\right\rangle+\left\langle d f_{\cdot, \widetilde{\mathcal{V}}}(X), \widetilde{\mathcal{W}}_{d f_{\mathcal{H}, \widetilde{\mathcal{H}}}\left(e_i\right)} d f_{\mathcal{H}, \widetilde{\mathcal{V}}}\left(e_i\right)\right\rangle \notag\\
    & =\left\langle d f_{\mathcal{H}, \widetilde{\mathcal{V}}}\left(e_i\right), \widetilde{\mathcal{T}}_{d f_{\mathcal{H}, \widetilde{\mathcal{H}}}\left(e_i\right)} d f_{, \widetilde{\mathcal{H}}}(X)\right\rangle+\left\langle d f_{\cdot, \widetilde{\mathcal{V}}}(X), \widetilde{\mathcal{W}}_{d f_{\mathcal{H}, \widetilde{\mathcal{H}}}\left(e_i\right)} d f_{\mathcal{H}, \widetilde{\mathcal{V}}}\left(e_i\right)\right\rangle  \notag\\
    & =\left\langle\widetilde{\mathcal{T}}^*_{d f_{\mathcal{H}, \widetilde{\mathcal{H}}}\left(e_i\right)} d f_{\mathcal{H}, \widetilde{\mathcal{V}}}\left(e_i\right), d f_{, \widetilde{\mathcal{H}}}(X)\right\rangle+\left\langle d f_{\cdot, \widetilde{\mathcal{V}}}(X), \widetilde{\mathcal{W}}_{d f_{\mathcal{H}, \widetilde{\mathcal{H}}}\left(e_i\right)} d f_{\mathcal{H}, \widetilde{\mathcal{V}}}\left(e_i\right)\right\rangle \notag\\
    & =\left\langle\widetilde{\mathcal{T}}^*_{d f_{\mathcal{H}, \widetilde{\mathcal{H}}}\left(e_i\right)} d f_{\mathcal{H}, \widetilde{\mathcal{V}}}\left(e_i\right)+\widetilde{\mathcal{W}}_{d f_{\mathcal{H}, \widetilde{\mathcal{H}}}\left(e_i\right)} d f_{\mathcal{H}, \widetilde{\mathcal{V}}}\left(e_i\right), d f(X)\right\rangle, \label{eq45}
\end{align}
    where the last equality in $\eqref{eq45}$ is due to the facts that $\widetilde{\mathcal{T}}^*$ and $\widetilde{\mathcal{W}}$ are $\widetilde{\mathcal{H}}$-valued and $\widetilde{\mathcal{V}}$-valued respectively. Using $\eqref{dfde}$ again, we obtain
    \begin{align}
         (I I I)=&\underbrace{\left\langle d f_{\mathcal{H}, \widetilde{\mathcal{H}}}\left(e_i\right), \widetilde{\nabla}_{e_i} d f(X)\right\rangle}_{(I I I)\text{-}1}-\underbrace{\left\langle d f_{\mathcal{H}, \widetilde{\mathcal{H}}}\left(e_i\right), \widetilde{\nabla}_{e_i} d f_{\mathcal{H}, \widetilde{\mathcal{V}}}(X)\right\rangle}_{(I I I)\text{-}2} \notag\\
        & -\underbrace{\left\langle d f_{\mathcal{H}, \widetilde{\mathcal{H}}}\left(e_i\right), \widetilde{\nabla}_{e_i} d f_{\mathcal{V}, \widetilde{\mathcal{V}}}(X)\right\rangle}_{(I I I)\text{-}3}-\underbrace{\left\langle d f_{\mathcal{H}, \widetilde{\mathcal{H}}}\left(e_i\right), \widetilde{\nabla}_{e_i} d f_{\mathcal{V}, \widetilde{\mathcal{H}}}(X)\right\rangle}_{(I I I)\text{-}4} \notag\\
        & =\underbrace{\left\langle d f_{\mathcal{H}, \widetilde{\mathcal{H}}}\left(e_i\right), \widetilde{\nabla}_{e_i} d f(X)\right\rangle}_{(I I I)\text{-}1}+\underbrace{\left\langle\widetilde{\nabla}_{e_i} d f_{\mathcal{H}, \widetilde{\mathcal{H}}}\left(e_i\right), d f_{\mathcal{H}, \widetilde{\mathcal{V}}}(X)\right\rangle}_{(I I I)\text{-}2} \notag\\
        & +\underbrace{\left\langle\widetilde{\nabla}_{e_i} d f_{\mathcal{H}, \widetilde{\mathcal{H}}}\left(e_i\right), d f_{\mathcal{V}, \widetilde{\mathcal{V}}}(X)\right\rangle}_{(I I I)\text{-}3}-\underbrace{\left\langle d f_{\mathcal{H}, \widetilde{\mathcal{H}}}\left(e_i\right), \widetilde{\nabla}_{e_i} d f_{\mathcal{V}, \widetilde{\mathcal{H}}}(X)\right\rangle}_{(I I I)\text{-}4}.
    \end{align}
Similarly, we have 
\begin{equation}\label{eq47}
    \begin{aligned}
        (I V) & =\underbrace{\left\langle d f_{\mathcal{H}, \widetilde{\mathcal{H}}}\left(e_i\right), d f_{\mathcal{H}, \widetilde{\mathcal{H}}}\left(\nabla_{e_i} X\right)\right\rangle}_{(I V)} \\
        & =\underbrace{\left\langle d f_{\mathcal{H}, \widetilde{\mathcal{H}}}\left(e_i\right), d f\left(\nabla_{e_i} X\right)\right\rangle}_{(I V)-1}-\underbrace{\left\langle d f_{\mathcal{H}, \widetilde{\mathcal{H}}}\left(e_i\right), d f_{\mathcal{V}, \widetilde{\mathcal{H}}}\left(\nabla_{e_i} X\right)\right\rangle}_{(I V)-2} .
    \end{aligned}
\end{equation}
Since $(II)\text{-}1 = (III)\text{-}1, (II)\text{-}3 = (IV )\text{-}1$, we obtain from $\eqref{divs}$-$\eqref{eq47}$ that
    \begin{align}
         \left(\operatorname{div} S_{\mathcal{H}, \widetilde{\mathcal{H}}}(f)\right)(X)=&-\underbrace{\left\langle\operatorname{tr}_g \beta_{\mathcal{H}, \widetilde{\mathcal{H}}}(f), d f(X)\right\rangle}_{(I)\text{-}1} \notag\\
        & +\underbrace{\left\langle\operatorname{tr}_g \beta_{\mathcal{H}, \widetilde{\mathcal{H}}}(f), d f_{\mathcal{V}, \widetilde{\mathcal{H}}}(X)+d f_{\mathcal{H}, \widetilde{\mathcal{V}}}(X)+d f_{\mathcal{V}, \widetilde{\mathcal{V}}}(X)\right\rangle}_{(I)\text{-}2} \notag\\
        & +\underbrace{\left\langle-d f_{\mathcal{H}, \widetilde{\mathcal{H}}}\left(\mathcal{T}_{e_i}^*(X)\right)+d f_{\mathcal{V}, \widetilde{\mathcal{H}}}\left(\mathcal{T}_X e_i\right)-d f_{\mathcal{V}, \widetilde{\mathcal{H}}}\left(\mathcal{W}_{e_i}(X)\right), d f_{\mathcal{H}, \widetilde{\mathcal{H}}}\left(e_i\right)\right\rangle}_{(I I)\text{-}2} \notag\\
        & -\underbrace{\left\langle\operatorname{tr}_g f^* \widetilde{\mathcal{T}}+t r_g f^* \widetilde{\mathcal{W}}, d f(X)\right\rangle}_{(I I)\text{-}4}-\underbrace{\left\langle\operatorname{tr}_g \beta_{\mathcal{H}, \widetilde{\mathcal{H}}}(f), d f_{\mathcal{H}, \widetilde{\mathcal{V}}}(X)+d f_{\mathcal{V}, \widetilde{\mathcal{V}}}(X)\right\rangle}_{(I I I)\text{-}2+(I I I)\text{-}3} \notag\\
        & +\underbrace{\left\langle d f_{\mathcal{H}, \widetilde{\mathcal{H}}}\left(e_i\right),\left(\widetilde{\nabla}_{e_i} d f_{\mathcal{V}, \widetilde{\mathcal{H}}}\right)(X)\right\rangle}_{(I I I)\text{-}4-(I V)\text{-}2} \notag\\
         =&-\underbrace{\left\langle\tau_{\mathcal{H}, \widetilde{\mathcal{H}}}(f), d f(X)\right\rangle}_{(I)\text{-}1+(I I)\text{-}4}+\underbrace{\left\langle\operatorname{tr}_g \beta_{\mathcal{H}, \widetilde{\mathcal{H}}}(f), d f_{\mathcal{V}, \widetilde{\mathcal{H}}}(X)\right\rangle}_{(I)\text{-}2-[(I I I)\text{-}2+(I I I)\text{-}3]}+\underbrace{\left\langle d f_{\mathcal{H}, \widetilde{\mathcal{H}}}\left(e_i\right),\left(\widetilde{\nabla}_{e_i} d f_{\mathcal{V}, \widetilde{\mathcal{H}}}\right)(X)\right\rangle}_{(I I I)\text{-}4-(I V)\text{-}2} \notag\\
        & +\underbrace{\left\langle-d f_{\mathcal{H}, \widetilde{\mathcal{H}}}\left(\mathcal{T}_{e_i}^*(X)\right)+d f_{\mathcal{V}, \widetilde{\mathcal{H}}}\left(\mathcal{T}_X e_i\right)-d f_{\mathcal{V}, \widetilde{\mathcal{H}}}\left(\mathcal{W}_{e_i}(X)\right), d f_{\mathcal{H}, \widetilde{\mathcal{H}}}\left(e_i\right)\right\rangle}_{(I I)\text{-}2} . 
    \end{align}
Note that for any $X\in \mathfrak{X}(M)$:
\begin{equation}
    \begin{aligned}
        \la df_{\mathcal{H},\widetilde{\mathcal{H}}}(\mathcal{T}^*_{e_i}X), df_{\mathcal{H},\widetilde{\mathcal{H}}}(e_i)   \ra =&  \la df_{\mathcal{H},\widetilde{\mathcal{H}}}(\la \mathcal{T}^*_{e_i}X, e_j\ra e_j), df_{\mathcal{H},\widetilde{\mathcal{H}}}(e_i)   \ra\\
        =& \la \mathcal{T}^*_{e_i}X, e_j\ra \la df_{\mathcal{H},\widetilde{\mathcal{H}}}( e_j), df_{\mathcal{H},\widetilde{\mathcal{H}}}(e_i)   \ra\\
        =& \la X, \mathcal{T}_{e_i} e_j\ra \la df_{\mathcal{H},\widetilde{\mathcal{H}}}( e_j), df_{\mathcal{H},\widetilde{\mathcal{H}}}(e_i)   \ra\\
        =&0,
    \end{aligned} 
\end{equation}
since $\mathcal{T}$ is skew-symmetric.
\end{proof}

\begin{corollary}\label{cor14}
    Suppose $f:\left(M^{m+k}, \mathcal{F}^k , g\right) \rightarrow(N^{n+l}, \widetilde{\mathcal{F}}^l , \widetilde{g})$ is a horizontally harmonic map between two Riemannian foliations and $X \in \mathfrak{X}_{\mathcal{H}}$. Then
    $$
(\operatorname{div} S_{\mathcal{H}, \widetilde{\mathcal{H}}}(f))(X)=2\langle d f_{\mathcal{V}, \widetilde{\mathcal{H}}}\left(\mathcal{T}_X e_i\right), d f_{\mathcal{H}, \widetilde{\mathcal{H}}}\left(e_i\right)\rangle .
    $$
    In particular, if $\mathcal{T} \equiv 0$, then $(\operatorname{div} S_{\mathcal{H}, \widetilde{\mathcal{H}}}(f))(X)=0$ for any $X \in \mathfrak{X}_{\mathcal{H}}$.
\end{corollary}

\begin{corollary}\label{cor15}
     Suppose $f:(M^{m+k}, \mathcal{F}^k , g) \rightarrow (N^{n+l}, \widetilde{\mathcal{F}}^l , \widetilde{g})$ is a transversally harmonic map between two Riemannian foliations. Then
    $$
    (\operatorname{div} S_{\mathcal{H}, \widetilde{\mathcal{H}}}(f))(X)=0
    $$
    for any $X \in \mathfrak{X}(M)$.
\end{corollary}

We say that $f:\left(M^{m+k}, \mathcal{F}^k , g\right) \rightarrow(N^{n+l}, \widetilde{\mathcal{F}}^l , \widetilde{g})$ satisfies the transverse conservation law with respect to a vector field $X \in \mathfrak{X}(M)$ if
\begin{equation}\label{eq48}
(	\operatorname{div} S_{\mathcal{H}, \widetilde{\mathcal{H}}} (f))(X)=0.
\end{equation}
If $\eqref{eq48}$ holds for any $X \in \mathfrak{X}(M)$, then $f$ is said to satisfy the transverse conservation law.

For any $X \in \mathfrak{X}(M)$, we denote by $\theta_X$ its dual one form, that is,
\begin{equation}
    \theta_X(Y)=g(X, Y), \quad Y \in T M .
\end{equation}
The covariant derivative of $\theta_X$ is a 2-tensor field $\nabla \theta_X$ defined by
\begin{equation}
    \left(\nabla \theta_X\right)(Y, Z)=\left(\nabla_Z \theta_X\right)(Y)=g\left(\nabla_Z X, Y\right), \quad Y, Z \in T M .
\end{equation}
In particular, if $X=\operatorname{grad} u$ for some smooth function $u$, then $\theta_X=d u$, and $\nabla \theta_X=\operatorname{Hess}(u)$.

Let $\Theta$ be a symmetric 2-tensor field on $M$. A direct computation yields (cf. \cite{dongVanishingTheoremsVector2011a}, \cite{bairdStressenergyTensorsLichnerowicz2008})
\begin{equation}\label{eq51}
    \begin{aligned}
        \operatorname{div}\left(i_X \Theta\right) & =\left\langle\Theta, \nabla \theta_X\right\rangle+(\operatorname{div} \Theta)(X) \\
        & =\frac{1}{2}\left\langle\Theta, L_X g\right\rangle+(\operatorname{div} \Theta)(X),
    \end{aligned}
\end{equation}
where $X \in \mathfrak{X}(M)$. Let $D$ be any bounded domain of $M$ with piecewise $C^1$ boundary. By applying Stokes' theorem to $\eqref{eq51}$, we obtain
\begin{equation}\label{eq:stk:52}
    \begin{aligned}
        \int_{\partial D} \Theta(X, \nu) d S_g & =\int_D\left[\left\langle\Theta, \nabla \theta_X\right\rangle+(\operatorname{div} \Theta)(X)\right] d V_g \\
        & =\int_D\left[\frac{1}{2}\left\langle\Theta, L_X g\right\rangle+(\operatorname{div} \Theta)(X)\right] d V_g.
    \end{aligned}
\end{equation}
\begin{remark}
Note that in \cite{chiangTransversallyBiharmonicMaps2008}, Chiang and Wolak employed a different stress-energy tensor (i.e., the transverse stress-energy tensor) associated with the foliated map $f$. Their stress-energy tensor was defined by the induced map $\bar{f}$ 
between the base manifolds and was used to investigate transversally biharmonic maps.
Additional applications of the transverse stress-energy tensor can also be found in \cite{jungTransversallyHarmonicMaps2011b} and \cite{fu2023}. 
\end{remark}

\section{Monotonicity formulas}\label{se:monoformula}

\subsection{Mixed conformal  Euclidean spaces}\label{Monomixed}
~

In this section, we consider  horizontally harmonic maps from some mixed conformally flat Euclidean spaces (cf. Example \ref{exam4}), and establish a monotonicity inequality similar to that in \cite{jin_liouville_nodate}.

Let $(\mathbb{R}^{m+k},\mathbb{R}^{k},g)$ be the Riemannian foliation given in Example \ref{exam4}, that is, 
$g=\phi(x)g_{can}^{h}+\eta (x,y)g_{can}^{v}$, where $g_{can}^{h}$ and $g_{can}^{v}$ are canonical metrics on $\mathbb{R} ^{m} $ and $\mathbb{R} ^{k} $ respectively. Also denote by $g_{can}$ the canonical Euclidean metric on $\mathbb{R}^{m+k}$. We adopt the identification $\rr^{m+k}\approx \rr^{m}\times \rr^{k}$ as a manifold with Euclidean coordinates $(x^{1},...x^{m},y^{1},...,y^{k})$, and assume $m>2$ and $k\geq 1$, unless otherwise stated. Then $(\rr^{m+k},\rr^{k},g)$
has a global orthonormal frame field $\{e_{i},e_{m+\alpha }\}_{1\leq i\leq m,%
\text{ }1\leq \alpha \leq k}=\{\phi ^{-1/2}\frac{\partial }{\partial x^{i}}%
,\eta ^{-1/2}\frac{\partial }{\partial y^{\alpha }}\}_{1\leq i\leq m,\text{ }%
1\leq \alpha \leq k}$, and the volume form $dV_{g}=\phi ^{m/2}\eta
^{k/2}dx\wedge dy$, where $dx\wedge dy=dx^{1}\wedge \cdots dx^{m}\wedge
dy^{1}\wedge \cdots \wedge dy^{k}$.

Let $|\cdot | $ be the standard Euclidean norm. Set
\begin{equation}\label{cy}
    D_{\rho, \delta}\left(x_0, y_0\right)=
    \left\{(x, y) \in \mathbb{R}^{m+k}: x \in \rr^m, y \in \rr^k,\left| x-x_0 \right| \leq \rho,\left| y-y_0 \right|  \leq \delta\right\} 
\end{equation}
for any $(x_0,y_0)\in \rr^{m+k} $.
Clearly,
$$
\partial D_{\rho, \delta}\left(x_0, y_0\right)=C_{\rho, \delta}^{(1)} \cup C_{\rho, \delta}^{(2)} ,
$$
where
$$
\begin{aligned}
    & C_{\rho, \delta}^{(1)}=\left\{(x, y) \in \rr^{m+k}: x \in \rr^m, y \in \rr^k,\left| x-x_0 \right|=\rho,\left| y-y_0 \right| \leq \delta\right\}, \\
    & C_{\rho, \delta}^{(2)}=\left\{(x, y) \in \rr^{m+k}: x \in \rr^m, y \in \rr^k,\left| x-x_0 \right| \leq \rho,\left| y-y_0 \right|=\delta\right\} .
\end{aligned}
$$
The  volume form  $dS_g$ on $C_{\rho, \delta}^{(1)}$ is given by
\[dS_g=\phi^{\frac{m-1}{2}}\eta^{\frac{k}{2}}dS_{g_{can}}.\]
Now we can derive the monotonicity inequality for horizontally harmonic maps from $(\rr^{m+k},\rr^k,g) $.
\begin{lemma}\label{lm:mixdd}
    Let $  u : (\mathbb{R}^{m+k}, \mathbb{R}^k, g) \rightarrow (N^{n+l}, \widetilde{\mathcal{F}}^{l}, \widetilde{g})$ be a  horizontally harmonic map. 
    Let $r:\mathbb{R}^{m+k}=\mathbb{R}^{m}\times \mathbb{R}^{k}\rightarrow \rr$ be the function defined by
    $r(x,y)=| x-x_{0}|$    for any $(x,y)\in 
    \mathbb{R}^{m+k}$ and a fixed point $x_{0}\in 
    \mathbb{R}^{m}$.
    Suppose $\phi  $ and $\eta  $ satisfy the following assumption $( A_1)$: there exist $\sigma >0$ and $
    R_0 \geq  0$ such that
$$ \left( \frac{m-2}{2} r \frac{\partial \log \phi }{\partial r}+\frac{k}{2} r \frac{\partial \log \eta }{\partial r}   \right) \geq \sigma -m+2, \quad \text{for } r > R_0,$$
where   $r\frac{\partial }{\partial r}=\frac{1}{2}\nabla ^{0}r^{2}$,
and $\nabla ^{0}$ denotes the Levi-Civita connection of $(
\mathbb{R}^{m+k},g_{can})$.
    Then
    \begin{equation*}
        R^{-\sigma }\left\{\int_{D_{R,\delta }(x_{0},y_{0})\setminus D_{R_0,\delta }(x_{0},y_{0})} |\dtu|^2 dV_g+\sigma ^{-1}H(R_0)\right\}
    \end{equation*}
    is an increasing function of $R$ for $R>R_0$,   where $H(R_0):= \int_{C_{R_0,\delta}^{(1)}} S_{\mathcal{H}, \widetilde{\mathcal{H}}}(u)\left(r\frac{\partial}{\partial r}, \nu\right) d S_g$.  
\end{lemma}

\begin{proof}
 Let $X=r\frac{\partial }{\partial r}$. Clearly, $X\in \Gamma (\mathcal{H})$.

Set $D_{R,\delta }:=D_{R,\delta }(x_{0},y_{0})$. Applying $\eqref{eq:stk:52} $, Corollary
\ref{cor14} to $f$ on $D_{R,\delta }\backslash D_{R_{0},\delta }$, and noting that
\[
S_{\mathcal{H},\widetilde{\mathcal{H}}}(f)(r\frac{\partial }{\partial r},\nu )=0\text{ \ on }%
\partial \left( D_{R,\delta }\backslash D_{R_{0},\delta }\right) \text{ }%
\backslash \text{ }\left( C_{R,\delta }^{(1)}\cup C_{R_{0},\delta
}^{(1)}\right) ,
\]%
we have  
\begin{equation}\label{eq:mixstk}
	\begin{aligned}
		&\int_{C_{R, \delta}^{(1)} } S_{\mathcal{H}, \widetilde{\mathcal{H}}}(u)\left(r\frac{\partial}{\partial r}, \nu\right) d S_g-\int_{C_{R_0, \delta}^{(1)}} S_{\mathcal{H}, \widetilde{\mathcal{H}}}(u)\left(r\frac{\partial}{\partial r}, \nu\right) d S_g\\
		=& \int_{D_{R, \delta}\setminus D_{R_0, \delta}} g(S_{\HH}(u),\frac{1}{2}L_{r \frac{\partial }{\partial r} }g) d V_g.
	\end{aligned} 
\end{equation}
By the definition of $g$, we derive that  
\begin{equation*}
	\begin{aligned}
		 L_{r \frac{\partial }{\partial r} }g=&  L_{r \frac{\partial }{\partial r} }(\phi g^h+\eta g^v)\\
		&= r \frac{\partial \phi }{\partial r}g^h +r \frac{\partial \eta }{\partial r}g^v+ (\phi \cdot L_{r \frac{\partial }{\partial r} }g^h+\eta \cdot L_{r \frac{\partial }{\partial r} }g^v)  \\
		&= r \frac{\partial \log \phi }{\partial r}(\phi g^h)+ r \frac{\partial \log \eta }{\partial r}(\phi g^v) + (\phi \cdot L_{r \frac{\partial }{\partial r} }g^h+\eta \cdot L_{r \frac{\partial }{\partial r} }g^v). 
	\end{aligned} 
\end{equation*}
Obviously,
\begin{equation*}
	( L_{r \frac{\partial }{\partial r} }g^h)(e_i,e_j)=( L_{r \frac{\partial }{\partial r} }g_{can})(e_i,e_j)=\hess_{g_{can}}(r^2)(e_i,e_j), \quad \forall i,j=1,\dots,m.
\end{equation*}
 Thus     
\begin{equation*}
		g \left(S_{\HH}(u), \phi \cdot L_{r \frac{\partial }{\partial r} }g^h\right)=\phi \cdot  g \left(S_{\HH}(f),  \hess_{g_{can}}(r^2)\right).
\end{equation*}
On the other hand,
\begin{equation*}
	\begin{aligned}
		g \left(S_{\HH}(u), \eta \cdot L_{r \frac{\partial }{\partial r} }g^v\right)=& \frac{\eta }{2}|\dtu|^2( L_{r \frac{\partial }{\partial r} }g^v)(e_{m+\alpha},e_{m+\alpha }) \\
		=& \frac{1 }{2}|\dtu|^2( L_{r \frac{\partial }{\partial r} }g^v)(\frac{\partial }{\partial y^{\alpha }} ,\frac{\partial }{\partial y^{\alpha }} ) \\
		=&0.
	\end{aligned} 
\end{equation*}
Therefore,
\begin{equation}\label{eq:mix1}
    \begin{aligned}
        g \left(S_{\HH}(u),  L_{r \frac{\partial }{\partial r} }g\right)=& r \frac{\partial \log \phi }{\partial r} \cdot g \left(S_{\HH}(u), \phi g^h\right)+r \frac{\partial \log \eta }{\partial r} \cdot g \left(S_{\HH}(u), \eta g^v\right) \\
       &+\phi \cdot  g \left(S_{\HH}(u),  \hess_{g_{can}}(r^2)\right).
    \end{aligned} 
\end{equation}
Note that 
\begin{equation}\label{eq:mix2}
	\begin{aligned}
		g \left(S_{\HH}(u), \phi g^h\right)=& \frac{m }{2} |\dtu|^2-|\dtu|^2=(m-2)\frac{|\dtu|^2}{2},\\
		g \left(S_{\HH}(u), \eta g^v\right)=& \frac{k }{2} |\dtu|^2
	\end{aligned} 
\end{equation}
and 
\begin{equation}\label{eq:mix3}
	\begin{aligned}
		\phi g \left(S_{\HH}(u),  \hess_{g_{can}}(r^2)\right)=& \frac{|\dtu|^2}{2}\Delta _{g_{can}}(r^2)\\
		&-\la \dtu(e_i),\dtu(e_j) \ra \hess_{g_{can}}(r^2)(\frac{\partial }{\partial x^i} ,\frac{\partial }{\partial x^j} )\\
		=& (m-2)|\dtu|^2.
	\end{aligned} 
\end{equation}
Putting $\eqref{eq:mix2},\eqref{eq:mix3} $ into $\eqref{eq:mix1} $, we get 
\begin{equation*}
	g \left(S_{\HH}(u),  L_{r \frac{\partial }{\partial r} }g\right)=\left[ (m-2)r \frac{\partial \log \phi }{\partial r}+kr \frac{\partial \log \eta }{\partial r}  +2m-4\right] \frac{|\dtu|^2}{2}.
\end{equation*}

On the other hand, by the coarea formula and $|\nabla r| =\phi ^{-\frac{1}{2}}$, we get  
\begin{equation*}
	\begin{aligned}
		\int_{C_{R,\delta}^{(1)}} S_{\mathcal{H}, \widetilde{\mathcal{H}}}(u)\left(r\frac{\partial}{\partial r}, \nu\right) d S_g=&\int_{C_{R,\delta}^{(1)}} S_{\mathcal{H}, \widetilde{\mathcal{H}}}(u)\left(r\frac{\partial}{\partial r}, \phi^{-\frac{1}{2}}\frac{\partial }{\partial r}  \right) d S_g\\
		=&\int_{C_{R,\delta}^{(1)}} r\frac{|\dtu|^2}{2} \phi ^{\frac{1}{2}}-r \phi ^{-\frac{1}{2}}\la \dtu(\frac{\partial }{\partial r} ), \dtu(\frac{\partial }{\partial r} )\ra  d S_g\\
		\leq & R\int_{C_{R,\delta}^{(1)}} \frac{|\dtu|^2}{2} \phi ^{\frac{1}{2}} dS_g\\
		=& R \frac{d}{dR} \int_{0}^{R}\left[ \frac{\int_{C_{t,\delta}^{(1)}} \frac{|\dtu|^2}{2}  dS_g}{|\nabla r|}\right]dt\\
		=& R \frac{d}{dR} \int_{D_{R,\delta }} \frac{|\dtu|^2}{2}dV_g.
	\end{aligned} 
\end{equation*}
Therefore, by $\eqref{eq:mixstk} $,
$$ R \frac{d}{dR} \int_{D_{R,\delta }} \frac{|\dtu|^2}{2}dt -H(R_{0})\geq  \frac{1}{2} \int_{D_{R,\delta }} \left( (m-2)r \frac{\partial \log \phi }{\partial r}+kr \frac{\partial \log \eta }{\partial r}   +2m-4\right)\frac{|\dtu|^2}{2} dV_g.$$
Using the assumption $( A_1)$,
we have 
$$ R \frac{d}{dR} \int_{D_{R,\delta }} \frac{|\dtu|^2}{2}dV_g -H(R_0)\geq  \sigma  \int_{D_{R,\delta }\setminus D_{R_0,\delta }} \frac{|\dtu|^2}{2} dV_g. $$  
It follows that 
$$ \frac{d}{dR} \frac{\int_{D_{R,\delta }\setminus D_{R_0,\delta }} |\dtu|^2 dV_g+\sigma ^{-1}H(R_0)}{R^{\sigma }} \geq 0, \quad \text{for }  R > R_0. $$
\end{proof}

\begin{remark}
    If the assumption $( A_1)$ holds for $R_0=0$, then $H=0$.    
\end{remark}

Taking $\phi \equiv 1, \eta \equiv 1 $ and $\sigma =m-2 $ in Lemma \ref{lm:mixdd}, we get  
\begin{theorem}\label{thm:eud}
     Let $f:(\mathbb{R}^{m+k}, \mathbb{R}^k, g_{c a n}) \rightarrow(N, \widetilde{\mathcal{F}}, \widetilde{g})$ be a horizontally harmonic map. Assume that $m>2$. Then
    $$
    \frac{\int_{D_{\rho_1, \delta}(x_0, y_0)}\left|d f_{\mathcal{H}, \widetilde{\mathcal{H}}}\right|^2 d V_g}{\rho_1^{m-2}} \leq \frac{\int_{D_{\rho_2, \delta}\left(x_0, y_0\right)}\left|d f_{\mathcal{H}, \widetilde{\mathcal{H}}}\right|^2 d V_g}{\rho_2^{m-2}}
    $$
    for any $\left(x_0, y_0\right) \in \mathbb{R}^{m+k}, \delta>0$ and $0<\rho_1 \leq \rho_2$.
\end{theorem}

\begin{remark}
    Note that  $ D_{\rho, \delta}\left(x_0, y_0\right)$ is a kind of ``cylinder'' defined by the distance to the leaf passing through $(x_0, y_0)$. We will use similar notation for general Riemannian foliations in the following.
\end{remark}

From Theorem \ref{thm:eud}, we obtain immediately the following lemma and theorem.
\begin{lemma}
 Let $f:(\mathbb{R}^{m+k}, \mathbb{R}^k, g_{can}) \rightarrow(N, \widetilde{\mathcal{F}}, \widetilde{g})$ be a horizontally harmonic map with $m>2$. If $f$ is not horizontally constant, then
    $$
    \int_{D_{\rho, \delta_0}\left(x_0, y_0\right)}\left|d f_{\mathcal{H}, \widetilde{\mathcal{H}}}\right|^2 d V_g \geq c\left(f, \delta_0\right) \rho^{m-2} \quad \text { as } \rho \rightarrow \infty \text {, }
    $$
    for some $\delta_0>0$, where $c\left(f, \delta_0\right)$ is a constant only depending on $f$ and $\delta_0$. In particular, the horizontal energy $E_{\mathcal{H} , \widetilde{\mathcal{H}}}(f)$ is infinite.
\end{lemma}

\begin{theorem}
    Let $f:(\mathbb{R}^{m+k}, \mathbb{R}^k, g_{c a n}) \rightarrow(N, \widetilde{\mathcal{F}}, \widetilde{g})$ be a horizontally harmonic map. Assume that $m>2$. If
    $$
    \int_{D_{\rho, \delta}\left(x_0, y_0\right)}\left|d f_{\mathcal{H}, \widetilde{\mathcal{H}}}\right|^2 d V_g=o\left(\rho^{m-2}\right) \quad \text { as } \rho \rightarrow \infty
    $$
    for any $\delta>0$, then $d f_{\mathcal{H}, \widetilde{\mathcal{H}}} \equiv 0$, that is, $f\left(\mathbb{R}^m \times\{p\}\right)$ is contained in a single leaf for any $p \in \mathbb{R}^k$.
\end{theorem}
\,

\subsection{The quotient space $K_m$ of the Heisenberg group }
~

In the previous subsection we considered the horizontally harmonic map from
the simplest model, i.e. the Euclidean space foliated by Euclidean subspaces, where the horizontal distribution
is integrable. In this subsection, let us consider the opposite case, i.e.,
the horizontally harmonic maps from the quotient space $K_{m}$ of the
Heisenberg group $H_{m}$ (cf. Example \ref{ex:sasa}), where $K_{m}$ is a simple model
of Riemannian foliation with non-integrable horizontal distribution.
Apparently, the Webster metric $g$ of $H_{m}$ descends in a natural way to a
metric on $K_{m}$ still denoted by $g$. 

From Example \ref{ex:sasa} (cf. also \cite{dragomirDifferentialGeometryAnalysis2006}), it is easy to verify that 
\begin{equation*}
X_j = \frac{\partial}{\partial x^j} + 2y^j \frac{\partial}{\partial t}, \quad Y_j = \frac{\partial}{\partial y^j} - 2x^j \frac{\partial}{\partial t}, \quad T =2 \frac{\partial}{\partial t}
\end{equation*}
 form a global orthonormal frame field on $K_{m}=%
\mathbb{C}
^{m}\times S^{1}$, where $t$ is the angle coordinate on $S^1$. The (CR) complex structure $J$ on $K_{m}$
is defined by $J(X_i)=Y_i, J(Y_i)=-X_i$ and $J(T)=0.$ Clearly $\mathcal{H}=\operatorname{span}\{X_1, \ldots, X_m, Y_1, \ldots, Y_m\}$ and $\mathcal{V}%
=\operatorname{span}\left\{T\right\}$. The Lie bracket relations of $\{X_{j},Y_{j},T\}$ are given by
\begin{equation*}
[Y_j, X_k] = 2\delta_{jk} T, \quad [X_j, X_k] = [Y_j, Y_k] = [X_j, T] = [Y_j, T] = 0.
\end{equation*}
Then Koszul's formula yields 
\begin{equation}\label{eq:h10}
    \begin{aligned}
        &\nabla  _{Y_j}X_i=-\nabla  _{X_i}Y_j=\delta _{ij}T, \quad \nabla  _{T}X_j=\nabla  _{X_j}T=Y_j,\\
        &\nabla  _{T}Y_j=\nabla  _{Y_j}T=-X_j,\quad  \nabla  _{T}T=\nabla  _{X_i}X_j=\nabla  _{Y_i}Y_j=0,
    \end{aligned} 
\end{equation}
where $\nabla $ denotes the Levi-Civita connection of $g$.

Next, for any $X=a^iX_i+b^jY_j \in \mathfrak{X}_{\mathcal{H}}$, there holds 
$$ \mathcal{T}_XX_i=\frac{1}{2}[X,X_i]^v=b^iT,\quad \mathcal{T}_XY_i=\frac{1}{2}[X,Y_i]^v=-a^iT.  $$ 
Let $f:\left(K_m, \mathcal{F}, g, J\right) \rightarrow(N^{n+l}, \widetilde{\mathcal{F}}^l, \widetilde{g})$ be a horizontally harmonic map from  $K_m$ to a Riemannian foliation $N$.
Therefore,  we have
\begin{equation}
    \begin{aligned}
        \left(\operatorname{div} S_{\mathcal{H}, \widetilde{\mathcal{H}}}(f)\right)(X)=&\left\langle 2 d f_{\mathcal{V}, \widetilde{\mathcal{H}}}\left(\mathcal{T}_X X_i\right), d f_{\mathcal{H}, \widetilde{\mathcal{H}}}\left(X_i\right)\right\rangle+\left\langle 2 d f_{\mathcal{V}, \widetilde{\mathcal{H}}}\left(\mathcal{T}_X Y_i\right), d f_{\mathcal{H}, \widetilde{\mathcal{H}}}\left(Y_i\right)\right\rangle\\
        =&2b^i \left\langle  d f_{\mathcal{V}, \widetilde{\mathcal{H}}}(T), d f_{\mathcal{H}, \widetilde{\mathcal{H}}}\left(X_i\right)\right\rangle -2a^i \left\langle  d f_{\mathcal{V}, \widetilde{\mathcal{H}}}(T), d f_{\mathcal{H}, \widetilde{\mathcal{H}}}\left(Y_i\right)\right\rangle \\
        =& 2 \left\langle  d f_{\mathcal{V}, \widetilde{\mathcal{H}}}(T), d f_{\mathcal{H}, \widetilde{\mathcal{H}}}\left(b^iX_i-a^jY_j\right)\right\rangle\\
        =& -2 \left\langle  d f_{\mathcal{V}, \widetilde{\mathcal{H}}}(T), d f_{\mathcal{H}, \widetilde{\mathcal{H}}}\left(JX\right)\right\rangle.
    \end{aligned} 
\end{equation}
Set
$$ D_{\rho}=\left\{(z, t) \in K_m: z \in \mathbb{C}^m, t \in S^1,\left|z\right| \leq \rho \right\}. $$
Then 
$$ \begin{aligned}
	& \partial D_{\rho } =\left\{(z, t) \in K_m: z \in \mathbb{C}^m, t \in S^1,|z|=\rho\right\} .
	\end{aligned} $$
Choose a basis $\{e_1,...,e_{2m-1}, e_{2m}= \frac{\partial }{\partial r} \}$ on the base manifold $\cc^m$, where $r(x,y)$ is the distance between $(x,y)$ and the origin $o$ of $\cc^m$. Let $E_1,...,E_{2m-1}, E_{2m} $ be their horizontal lifts. Clearly $\{E_1,...,E_{2m-1}, E_{2m},T\}$ constitutes a basis for $K_m$.         

Now take $X=(r\circ\pi)\cdot E_{2m} \in  \mathfrak{X}_{\mathcal{H}}$, that is, 
$$ X=\sqrt{|x|^2+|y|^2}\left(\frac{x^i}{r} X_i+\frac{y^i}{r} Y_i\right)=x^i X_i+y^j Y_j. $$
It follows that 
$$ \mathcal{T}_{X}X_i=y^iT,\quad  \mathcal{T}_{X}Y_i=-x^iT.$$
In addition,  by $\eqref{ttensor}$  and $\eqref{eq:h10}$,
$$
\begin{aligned}
\mathcal{T}_X E_a & =\left\langle\nabla_X E_a, T\right\rangle T \\
& =-\left\langle E_a, \nabla_X T\right\rangle T \\
& =-\left\langle E_a, x^i Y_i-y^j X_j\right\rangle T \\
& =-\left\langle E_a, J X\right\rangle T,
\end{aligned}
$$
where $1 \leq a \leq 2m$.
Thus, due to Corollary \ref{cor14},
\begin{equation}\label{eq:h95}
	(\di \st)(X)=-2r \la df_{\mathcal{V},\widetilde{\mathcal{H}}}(T),  \dt (J E_{2m}) \ra.  
\end{equation}

Noting that $\nu:=E_{2m}$ is the unit normal vector of $\partial D_{\rho}$, we have
\begin{equation}\label{eq:h96}
	S_{\mathcal{H}, \widetilde{\mathcal{H}}}(f)\left(X, \nu\right)=\rho e_{\mathcal{H}, \widetilde{\mathcal{H}}}(f)-\rho\left|d f_{\mathcal{H}, \widetilde{\mathcal{H}}}\left(E_{2m}\right)\right|^2 \quad \text{on} \quad \partial D_{\rho},
\end{equation}
where $e_{\mathcal{H}, \widetilde{\mathcal{H}}}(f)=\frac{1}{2}\left|d f_{\mathcal{H}, \widetilde{\mathcal{H}}}\right|^2$ denotes the horizontal energy density of $f$.

A direct computation shows that
$$
\begin{aligned}
\nabla_{E_b} X & =r \nabla_{E_b} E_{2m}, \quad 1 \leq b \leq 2m-1; \\
\nabla_{E_{2m}} X & =E_{2m}+r \nabla_{E_{2m}} E_{2m}; \\
\nabla_T X & =r \nabla_T E_{2m} .
\end{aligned}
$$
Noting that the foliation is minimal, we have
$$
\begin{aligned}
&\sum_{b=1}^{2m-1} \left\langle\nabla_{E_b} X, E_b\right\rangle+\left\langle\nabla_{E_{2m}} X, E_{2m}\right\rangle+\left\langle\nabla_T X, T\right\rangle \\
 =&r \sum_{b=1}^{2m-1} \left\langle\nabla_{E_b} E_{2m}, E_b\right\rangle+1+0 \\
 =&\sum_{b=1}^{2m-1} \operatorname{Hess}\left(\frac{r^2}{2}\right)\left(e_b, e_b\right)+1 \\
 =&2 m.
\end{aligned}
$$
Similarly, 
$$ \begin{aligned}
	\left\langle d f_{\mathcal{H}, \widetilde{\mathcal{H}}} \odot d f_{\mathcal{H}, \widetilde{\mathcal{H}}}, \nabla \theta_ X\right\rangle  =&\sum_{a,b=1}^{2m-1} \left\langle d f_{\mathcal{H}, \widetilde{\mathcal{H}}}\left(E_a\right), d f_{\mathcal{H}, \widetilde{\mathcal{H}}}\left(E_b\right)\right\rangle\left\langle\nabla_{E_a} X, E_b\right\rangle\\
 &+ \left\langle d f_{\mathcal{H}, \widetilde{\mathcal{H}}}\left(E_{2m}\right), d f_{\mathcal{H}, \widetilde{\mathcal{H}}}\left(E_{2m}\right)\right\rangle\left\langle\nabla_{E_{2m}} X, E_{2m}\right\rangle \\
	 =&\sum_{\widetilde{k}=1}^{n} \sum_{a, b=1}^{2 m-1} \operatorname{Hess}\left(\frac{r^2}{2}\right)\left(f_a^{\widetilde{k}} e_a, f_b^{\widetilde{k}} e_b\right)+\left|d f_{\mathcal{H}, \widetilde{\mathcal{H}}}\left(E_{2m}\right)\right|^2 \\
	=&\left|d f_{\mathcal{H}, \widetilde{\mathcal{H}}}\right|^2 .
	\end{aligned} $$
	Thus, 
	\begin{equation}\label{eq:h97}
		\left\langle S_{\mathcal{H}, \widetilde{\mathcal{H}}}(f), \nabla \theta_X\right\rangle=(2 m-2) e_{\mathcal{H}, \widetilde{\mathcal{H}}}(f) . 
	\end{equation}
Now we are ready to derive the monotonicity formula for $f$.
\begin{lemma}\label{le:hemo}
	Let $f:\left(K_m, \mathcal{F}, g, J\right) \rightarrow(N, \widetilde{\mathcal{F}}, \widetilde{g})$ be a horizontally harmonic map from  $K_m$ to a Riemannian foliation $N$. 
	Suppose that there are constants $\delta $ and $\rho _{0}>0$ such that    
	\begin{equation}\label{eq:condH}
	2\rho \left( \int_{D_{\rho }}|df_{\mathcal{V},\widetilde{\mathcal{H}}%
	}|^{2}dV_g \right) ^{1/2}\leq \delta \left( \int_{D_{\rho }}|df_{\mathcal{H},%
	\widetilde{\mathcal{H}}}|^{2}dV_g \right) ^{1/2}\text{ \ \ \ for \ }\rho \geq
	\rho _{0}  
	\end{equation}%
	and 
	\[
	2m-2-2\delta >0.
	\]%
    Then 
   $$
   \sigma^{2-2 m+2 \delta } \int_{D_{\sigma   }} e_{\mathcal{H}, \widetilde{\mathcal{H}}} (f) d V_g \leq \rho^{2-2 m+2 \delta} \int_{D_{\rho  }} e_{\mathcal{H}, \widetilde{\mathcal{H}}} (f) d V_g
   $$
   holds for any $\rho_0<\sigma \leq \rho$. 
\end{lemma}
\begin{proof}
   Putting $\eqref{eq:h95}, \eqref{eq:h96}$ and $\eqref{eq:h97}$ into the following  Stokes' formula 
$$ \int_{\partial D_{\rho  }} \st(X,\nu )dS_g =\int_{D_{\rho  }}[\la \st,\nabla \theta _X \ra+(\di \st)(X)  ]dV_g,$$ 
we conclude that
\begin{equation*}
   \begin{aligned}
   &\rho \int_{\partial D_{\rho  }  } \left[ e_{\mathcal{H}, \widetilde{\mathcal{H}}}(f)-\left|d f_{\mathcal{H}, \widetilde{\mathcal{H}}}\left(E_{2m}\right)\right|^2 \right]d S_g \\
   = & (2m-2)\int_{D_{\rho  }}  e_{\mathcal{H}, \widetilde{\mathcal{H}}}(f)d V_g -\int_{D_{\rho  }} 2 r\left\langle d f_{\mathcal{V}, \widetilde{\mathcal{H}}}(T),  d f_{\mathcal{H}, \widetilde{\mathcal{H}}}\left(J E_{2m}\right)\right\rangle d V_g.
   \end{aligned}
\end{equation*}
By H\"older inequality and $\eqref{eq:condH}$, we get
$$ \begin{aligned}
	\int_{D_{\rho }}2r\langle df_{\mathcal{V},\widetilde{\mathcal{H}}}(T),df_{%
\mathcal{H},\widetilde{\mathcal{H}}}(JE_{2m})\rangle dV_g\leq &2\rho \left(
\int_{D_{\rho }}|df_{\mathcal{V},\widetilde{\mathcal{H}}}|^{2}dV_g\right)
^{1/2}\left( \int _{D_{\rho }}|df_{\mathcal{H},\widetilde{\mathcal{H}}}|^{2}dV_g\right)
^{1/2}\\
	\leq & \delta \int_{D_{\rho }}| \dt |^2 dV_g,
\end{aligned}  $$
according to $ \eqref{eq:condH} $.
Therefore, 
$$ \begin{aligned}
	\rho \int_{\partial D_{\rho  }} e_{\mathcal{H}, \widetilde{\mathcal{H}}}(f) dS_g \geq & (2 m-2-2 \delta )\int_{D_{\rho  }} e_{\mathcal{H}, \widetilde{\mathcal{H}}}(f) d V_g ,
	\end{aligned} $$
	which implies 
	$$ \frac{d}{d \rho}\left(\rho^{-2 m+2+2 \delta } \int_{D_{\rho  }} e_{\mathcal{H}, \widetilde{\mathcal{H}}}(f) d V_g\right) \geq 0  \qquad \text{for} \quad \rho \geq \rho _0. $$
\end{proof}
\begin{remark}
    If
    $$2\rho \left( \int_{D_{\rho }}|df_{\mathcal{V},\widetilde{\mathcal{H}}%
	}|^{2}dV_g\right) ^{1/2}\leq C(\rho) \left( \int_{D_{\rho }}|df_{\mathcal{H},%
	\widetilde{\mathcal{H}}}|^{2}dV_g\right) ^{1/2}  $$
 with $C(\rho) \to 0$ as $ \rho \to \infty$, then there are constants $\delta$ and $\rho_0>0$ such that 
 $$C(\rho) < \delta \text{ \ \ \ for \ }\rho \geq
	\rho _{0} , $$
 and 
 $$2m-2-2\delta>0. $$
\end{remark}

The following theorem holds as an immediate result.
\begin{theorem}
	Assume that $f$ satisfies the same conditions as in Lemma \ref{le:hemo} and that 
	$$ \int_{D_{\rho }}\left| \dt \right|^2 dV_g=o(\rho ^{2m-2-2 \delta  }) \quad \text{as }  \, \rho \to +\infty.  $$ 
	Then $f$ is horizontally constant.  
	\end{theorem}

    \subsection{Transversally harmonic maps under curvature conditions }\label{transverHarm}

    \begin{lemma}\label{lm:gwtype}
        Let $(M^{m+k}, g)$ be a complete Riemannian manifold with a pole $x_0$ and let $r$ be the distance function relative to $x_0$. Suppose the radial curvature $K_r$ of $M$ satisfies one of the following:
        \begin{enumerate}[(i)]
            \item $-\alpha ^2\leq  K_r \leq -\beta ^2,  $ where $\alpha ,\beta >0 $ and $(m+k-1)\beta -2\alpha \geq 0; $
            \item   $K_r=0$ and $m+k-2 >0$ ;
            \item   $-\frac{A}{\left(1+r^2\right)^{1+\epsilon}} \leq K_r \leq \frac{B}{\left(1+r^2\right)^{1+\epsilon}}$ with $\epsilon>0, A \geq 0,0 \leq B<2 \epsilon$ and $m+k-(m+k-1)\frac{B}{2\epsilon}-2e^{\frac{A}{2 \epsilon }}>0.$
        \end{enumerate}  
        Suppose $(M^{m+k},\mathcal{F}^{k},g)$ and $(N,\widetilde{\mathcal{F}},%
    \widetilde{g})$ are Riemannian foliations and $f$ is a transversally
    harmonic map from $M$ to $N$. Then for any $0<\rho _1 \leq \rho _2 $, we have 
        \begin{equation}\label{eq:mono}
            \frac{1}{\rho_1^\lambda} \int_{B_{\rho _1}(x_0)} \frac{|\dt|^2}{2} dV_g \leq \frac{1}{\rho_2^\lambda} \int_{B_{\rho _2}(x_0)} \frac{|\dt|^2}{2} dV_g ,
        \end{equation}  
        where 
        \begin{equation*}
            \lambda = \begin{cases}
            m+k-\frac{2 \alpha }{\beta }	&\text{if}\quad K_r \,\text{satisfies (i)} \\
                m+k-2&\text{if} \quad K_r \, \text{satisfies (ii)}\\
                m+k-(m+k-1)\frac{B}{2\epsilon}-2e^{\frac{A}{2\epsilon}} &\text{if} \quad K_r \, \text{satisfies (iii).}
           \end{cases}   
        \end{equation*}
    \end{lemma}
    \begin{proof}
    Since $f$ is transversally harmonic, we have from Corollary \ref{cor15} that 
    \[
    \di S_{\mathcal{H},\widetilde{\mathcal{H}}}(f)=0.
    \]%
    As we have already mentioned in \S \ref{se:streeenergy}, $S_{\mathcal{H},\widetilde{\mathcal{H}}%
    }(f)\ $is the stress-energy tensor $S_{\omega }$, where $\omega =df_{%
    \mathcal{H},\widetilde{\mathcal{H}}}$ is the $f^{-1}\widetilde{\mathcal{H}}$%
    -valued $1$-form on $M$. Then Lemma \ref{lm:gwtype} follows immediately from Theorem
    4.1 in \cite{dongVanishingTheoremsVector2011a}.

    For the reader's convenience, we deduce the result of case (i) as follows.
    Assume that $-\alpha ^{2}\leq K_{r}\leq -\beta ^{2}$ with $\alpha ,\beta >0$
    and $(m+k-1)\beta -2\alpha $ $\geq 0$. It is known that (cf. \cite{greeneFunctionTheoryManifolds1979a}, \cite{dongVanishingTheoremsVector2011a})%
    \[
        \beta \operatorname{coth}(\beta r)[g-d r \otimes d r] \leq \operatorname{Hess}_g(r) \leq \alpha \operatorname{coth}(\alpha r)[g-d r \otimes d r].
    \]%
    Taking $X=r\frac{\partial }{\partial r}\in \mathfrak{X}(M)$, we derive that%
    $$ \begin{aligned}
        \la \st,\nabla \theta _X \ra \geq & \frac{\left| \dt \right|^2 }{2}( 1+(m+k-1)\beta r \operatorname{coth}(\beta r)-2\alpha r \operatorname{coth}(\alpha r)) \\
        =&  \frac{\left| \dt \right|^2 }{2}\left( 1+\beta r \operatorname{coth}(\beta r)[(m+k-1)-\frac{2\alpha r \operatorname{coth}(\alpha r)}{\beta r \operatorname{coth}(\beta r)}]\right)\\
        > & \frac{\left| \dt \right|^2 }{2}( 1+(m+k-1)-2\frac{\alpha }{\beta })\\
        =& \lambda \frac{\left| \dt \right|^2 }{2},
    \end{aligned}  $$
    where we have used the fact that
     $\beta r \coth (\beta r) $ is an increasing function which approaches to $1$ as $r \to 0$,  and $ \frac{\coth(\alpha r)}{\coth(\beta r)}<1$ (cf.  Lemma 4.2 in \cite{dongVanishingTheoremsVector2011a}).
     Applying $\eqref{eq:stk:52}$  to $\Theta =S_{\mathcal{H},%
     \widetilde{\mathcal{H}}}(f)$ and using the estimate $S_{\mathcal{H},%
     \widetilde{\mathcal{H}}}(f)(X,\frac{\partial }{\partial r})\leq \frac{\rho }{%
     2}|df_{\mathcal{H},\widetilde{\mathcal{H}}}|^{2}$ on $\partial B_{\rho }(x_{0})$, we get  
    $$ \rho \int_{\partial B_\rho (x_0)} \frac{\left| \dt \right|^2 }{2} dS_g \geq \lambda \int_{B_\rho (x_0)} \frac{\left| \dt \right|^2 }{2} dV_g. $$
    Then the co-area formula implies that 
    \begin{equation}\label{eq:pr65}
        \frac{ \frac{d}{d \rho } \int_{B_\rho (x_0)} \frac{\left| \dt \right|^2 }{2} dV_g}{ \int_{B_\rho (x_0)} \frac{\left| \dt \right|^2 }{2} dV_g} \geq \frac{\lambda }{\rho }.
    \end{equation}
    Now the monotonicity inequality $\eqref{eq:mono}$ follows by integrating  $\eqref{eq:pr65}$  on $[\rho _1,\rho _2]$.  
    \end{proof}

\section{Jin-type results}\label{se:Jintypethm}

\subsection{Mixed conformal  Euclidean spaces}\label{se:Mixedconformal}
~

In this subsection, based on the discussions in \S  \ref{Monomixed}, we establish a Jin
type theorem for a horizontally harmonic map $u:(
\mathbb{R}
^{m+k},
\mathbb{R}
^{k},g)\rightarrow (N^{n+l},\widetilde{F}^{l},\widetilde{g})$ from the mixed
conformal Euclidean space by further assuming that $u(x,y)\rightarrow q$ as $
(x,y)\rightarrow \infty $. Let $(U,\varphi ;\widetilde{x}^{\widetilde{A}})$
be a foliated chart around $q$ with $\varphi (U)=D_{\rho }^{n+l}$ as
described in \S \ref{se:Riemfoli}, which induces a Riemannian submersion $\pi :(U,\widetilde{
g})\rightarrow (\widetilde{U},\bar{h})$, and a coordinate chart $(
\widetilde{U},\widetilde{\varphi };\widetilde{x}^{\widetilde{i}})$ around $
\overline{q}=\pi (q)$ with $\widetilde{\varphi } \circ \pi 
=\Pr_{n}\circ \varphi $ and $\widetilde{\varphi }(\widetilde{U})=D_{\rho }^{n}$, 
such that $\widetilde{\varphi }(q)=o\in \rr^{n}$:
\[
\begin{tikzcd}
    U \arrow[r, "\varphi"] \arrow[d, "\pi"'] & D_{\rho }^{n+l} \arrow[d, "\Pr_n"] \\
    \widetilde{U} \arrow[r, "\widetilde{\varphi }"'] & D_{\rho }^{n}
\end{tikzcd}
\]

The asymptotic condition implies that there exists a sufficiently large $
R_{1}>0$ such that $u(x,y)\in U$ for $| (x,y)| >R_{1}$. Define a
smooth map $v:
\mathbb{R}
^{m+k}$ $\backslash $ $B_{R_{1}}\rightarrow \widetilde{U}$ by setting $
v:=\pi \circ u$. Clearly $\widetilde{\varphi }(v(x,y))\rightarrow o\in \rr^{n}$
as $(x,y)\rightarrow \infty $. The horizontal energy density of $u$ on $
\mathbb{R}
^{m+k}$ $\backslash $ $B_{R_{1}}$ is given by%
\[
e_{\mathcal{H},\widetilde{\mathcal{H}}}(u)=e_{\mathcal{H},\cdot }(v)=\frac{1
}{2}\phi ^{-1}\bar h_{\widetilde{i}\widetilde{j}}(v)\frac{\partial v^{\widetilde{i
}}}{\partial x_{i}}\frac{\partial v^{\widetilde{j}}}{\partial x_{i}}
\]
where $\bar{h}_{\widetilde{i}\widetilde{j}}= \bar{h}(\frac{\partial }{\partial \widetilde{x}^{\widetilde{i}}},\frac{\partial }{\partial \widetilde{x}^{\widetilde{j}}}  ) $.
From \S \ref{Monomixed}, we have   
\[
dV_{g}=V(x,y)dx\wedge dy,
\]
where $V(x,y):=\frac{1}{2}\eta ^{\frac{1}{2}}(x,y)\phi (x)^{\frac{m-2}{2}}$.
Then we can establish the following Jin-type theorem.

\begin{theorem}\label{thjin}
    Let $u : (\mathbb{R}^{m+k}, \mathbb{R}^k, g) \rightarrow (N^{n+l}, \widetilde{\mathcal{F}}^l, \widetilde{g})$ be a horizontally harmonic map with $m > 2, k\geq 1$ and $u(x,y)\to q \in N$ as $(x,y)\to \infty$. Suppose the condition $( A_1)$ holds for $R_0=0$ and  the following condition $( A_2)$ holds: 
    $$ V(x,y)\leq C_0(y)|x-x_0|^{\sigma-(m-2)}, $$
    where $\sigma $ is the constant in the condition $(A_1)$   and $C_0(y) $ is a smooth, positive function of $y$.
    Then  $u$ is horizontally constant. 
\end{theorem}

\begin{proof}   
	Let us suppose $du_{\mathcal{H},\widetilde{\mathcal{H}}}\neq 0$ in
order to derive a contradiction under the conditions of the theorem. In the following, we  denote $D_{R,\delta }:=D_{R,\delta }(x_{0},y_{0})$, and for simplicity, assume $(x_0,y_0)=(0,0)$.

Since $u(x,y)\rightarrow q$ and $v(x,y)\rightarrow \overline{q}$ as $|
(x,y)| \rightarrow \infty $, there is a sufficiently large $R_{1}>0$ such
that if $| (x,y)| >R_{1}$, then $u(x,y)\in U$, $v(x,y)\in \widetilde{U}
$ and 
\begin{equation}\label{eq71}
    \left( \frac{\partial \bar h_{{\widetilde{i}}{\widetilde{j}}}(v)}{\partial \widetilde{x}^{\widetilde{k}}} v^{\widetilde{k}} + 2 \bar h_{{\widetilde{i}}{\widetilde{j}}}(v) \right) \geq \left( h_{\widetilde{i}\widetilde{j}}(v)\right)   
\end{equation}
as matrices. For any $\widetilde{w}\in C_{0}^{2}\left( D_{\infty ,\delta
}\backslash D_{R_{1},\delta },\varphi  (U)\right) $ and sufficiently small $t$, we consider the variation $u+t\widetilde{w}:\mathbb{R}^{m+k}\rightarrow N$ defined as follows: 
\begin{equation}
    (u+t\widetilde{w})(x,y)=
    \begin{cases}
        u(x) &\text{if } (x,y)\in D_{R_{1},\delta }\\
        \varphi ^{-1}[\varphi \left( u(x,y)\right)+t\widetilde{w}(x,y)], &\text{if } (x,y)\in \mathbb{R}^{m+k}\backslash D_{R_{1},\delta }\\
   \end{cases} 
\end{equation}
Set $w=\Pr_n \circ \widetilde{w}\in C_{0}^{2}( D_{\infty ,\delta
}\backslash D_{R_{1},\delta }, \widetilde{\varphi } (\widetilde{U})) $. 
	It follows from the horizontal harmonicity of $u$ that:
    \[
    \left. \frac{d}{dt} \right|_{t=0} E_{\mathcal{H}, \widetilde{\mathcal{H}}}(u + t \widetilde{w}) = 0,
    \]
    that is,
    \begin{equation}\label{eq70}
        \int_{D_{\infty, \delta} \setminus D_{R_1, \delta}}  \left( 2\bar{h}_{{\widetilde{i}}{\widetilde{j}}}(v) \frac{\partial v^{\widetilde{i}}}{\partial x_i} \frac{\partial w^{\widetilde{j}}}{\partial x_i} + \frac{\partial \bar{h}_{{\widetilde{i}}{\widetilde{j}}}(v)}{\partial \widetilde{x}^{\widetilde{k}}} w^{\widetilde{k}} \frac{\partial v^{\widetilde{i}}}{\partial x_i} \frac{\partial v^{\widetilde{j}}}{\partial x_i} \right)  V(x,y)dx\wedge dy = 0.
    \end{equation}

    Choose \( \widetilde{w}(x,y) = \psi(|x|)\xi(|y|)\varphi (u(x,y)) \) in $\eqref{eq70} $  for  \( \psi(t) \in C_0^\infty(R_1, \infty) \) and $\xi(s) \in C_0^\infty(-\delta , \delta)$ with $\xi \equiv 1$ in $(-\delta /2,\delta/2)$. Since $ \widetilde{\varphi } \circ \pi 
    =\Pr_{n}\circ \varphi$,   $w(x,y) = \psi(|x|)\xi(|y|) \widetilde{\varphi } (v(x,y)) $.
    We have
\begin{equation}\label{eq50}
        \begin{aligned}
            &\int_{D_{\infty, \delta} \setminus D_{R_1, \delta}} \left( 2\bar{h}_{{\widetilde{i}}{\widetilde{j}}}(v) + \frac{\partial \bar{h}_{{\widetilde{i}}{\widetilde{j}}}(v)}{\partial \widetilde{x}^{\widetilde{k}}} v^{\widetilde{k}} \right) \frac{\partial v^{\widetilde{i}}}{\partial x_i} \frac{\partial v^{\widetilde{j}}}{\partial x_i} \psi(|x|) \xi(|y|)V(x,y) dx\wedge dy\\
            =& - \int_{D_{\infty, \delta} \setminus D_{R_1, \delta}} 2\bar{h}_{{\widetilde{i}}{\widetilde{j}}}(v) \frac{\partial v^{\widetilde{i}}}{\partial x_i} v^{\widetilde{j}} \frac{\partial \psi(|x|)}{\partial x_i} \xi(|y|)V(x,y) dx\wedge dy.
        \end{aligned}
    \end{equation}
    
    For $0<\varepsilon \leq 1$, define
    $$
    \zeta _{\varepsilon}(t)= \begin{cases}1 & t \leq 1; \\ 1+\frac{1-t}{\varepsilon} & 1<t<1+\varepsilon ; \\ 0 & t \geq 1+\varepsilon .\end{cases}
    $$
    In $\eqref{eq50}$, choose
    $$
    \psi(|x|)=\zeta_{\varepsilon}\left(\frac{|x|}{R}\right)\left(1-\zeta_1\left(\frac{|x|}{R_1}\right)\right) 
    $$
    for some $R> R_2 $,  where $R_2:=2R_1$.
    Notice that
    $$
    \frac{\partial \zeta_{\varepsilon}\left(\frac{|x|}{R}\right)}{\partial x_i}=-\frac{1}{R \varepsilon} \frac{x_i}{|x|} \quad  \text { for } R<|x|<R(1+\varepsilon),
    $$
    we get 
    \begin{equation}\label{eq52}
        \begin{aligned}
            & \int_{D_{R, \delta} \backslash D_{R_2, \delta}}\left(2 \bar{h}_{{\widetilde{i}} {\widetilde{j}}}(v)+\frac{\partial \bar{h}_{{\widetilde{i}} {\widetilde{j}}}(v)}{\partial \widetilde{x}^{\widetilde{k}}} v^{\widetilde{k}}\right) \frac{\partial v^{\widetilde{i}}}{\partial x_i} \frac{\partial v^{\widetilde{j}}}{\partial x_i}  \xi(|y|)V(x,y)dx\wedge dy+D\left(R_1\right) \\
            &+\int_{D_{\infty, \delta} \setminus D_{R, \delta}} \left( 2\bar{h}_{{\widetilde{i}}{\widetilde{j}}}(v) + \frac{\partial \bar{h}_{{\widetilde{i}}{\widetilde{j}}}(v)}{\partial \widetilde{x}^{\widetilde{k}}} v^{\widetilde{k}} \right) \frac{\partial v^{\widetilde{i}}}{\partial x_i} \frac{\partial v^{\widetilde{j}}}{\partial x_i} \zeta_{\varepsilon}\left(\frac{|x|}{R}\right)\xi(|y|)V(x,y)dx\wedge dy\\
            =&\frac{1}{R \varepsilon}\int_{D_{R(1+\varepsilon), \delta} \backslash D_{R, \delta}} 2 \bar{h}_{{\widetilde{i}} {\widetilde{j}}}(v) \frac{\partial v^{\widetilde{i}}}{\partial x_i}  v^{\widetilde{j}}   \frac{x_i}{|x|}\xi(|y|)V(x,y)dx\wedge dy,
        \end{aligned}
    \end{equation}
    where
    $$
    \begin{aligned}
        D\left(R_1\right)= & \int_{D_{R_2, \delta} \backslash D_{R_1, \delta}}\left(2 \bar{h}_{{\widetilde{i}} {\widetilde{j}}}(v)+\frac{\partial \bar{h}_{{\widetilde{i}} {\widetilde{j}}}(v)}{\partial \widetilde{x}^{\widetilde{k}}} v^{\widetilde{k}}\right) \frac{\partial v^{\widetilde{i}}}{\partial x_i} \frac{\partial v^{\widetilde{j}}}{\partial x_i}\left(1-\zeta_1\left(\frac{|x|}{R_1}\right)\right)  \xi(|y|)V(x,y)dx\wedge dy\\
        & -\int_{D_{R_2, \delta} \backslash D_{R_1, \delta}} 2 \bar{h}_{{\widetilde{i}} {\widetilde{j}}}(v) \frac{\partial v^{\widetilde{i}}}{\partial x_i} v^{\widetilde{j}} \frac{\partial \zeta_1\left(\frac{|x|}{R_1}\right)}{\partial x_i} \xi(|y|)V(x,y)dx\wedge dy.
    \end{aligned}
    $$
    Now, letting $\varepsilon \to 0$ in $\eqref{eq52}$, we have 
    \begin{equation}\label{eq53}
        \begin{aligned}
            & \int_{D_{R, \delta} \backslash D_{R_2, \delta}}\left(2 \bar{h}_{{\widetilde{i}} {\widetilde{j}}}(v)+\frac{\partial \bar{h}_{{\widetilde{i}} {\widetilde{j}}}(v)}{\partial \widetilde{x}^{\widetilde{k}}} v^{\widetilde{k}}\right) \frac{\partial v^{\widetilde{i}}}{\partial x_i} \frac{\partial v^{\widetilde{j}}}{\partial x_i}  \xi(|y|)V(x,y)dx\wedge dy+D\left(R_1\right) \\
            & =\int_{C^{(1)}_{R,\delta}} 2 \bar{h}_{{\widetilde{i}} {\widetilde{j}}}(v) \frac{\partial v^{\widetilde{i}}}{\partial x_i}  v^{\widetilde{j}}   \nu^i\xi(|y|) V(x,y)d S_{g_{can}} ,
        \end{aligned}
    \end{equation}
    where $\nu ^i=\frac{x_i}{|x|}$. Thus,   
    $\nu=\nu ^i \frac{\partial }{\partial x_i} $ is the unit outer normal vector field along $C^{(1)}_{R,\delta}$ relative to $g_{can}$.
    
    Set
    \begin{equation*}
        Z(R, \delta)=\int_{D_{R, \delta} \backslash D_{R_2, \delta}} \bar{h}_{{\widetilde{i}} {\widetilde{j}}}(v) \frac{\partial v^{\widetilde{i}}}{\partial x_i} \frac{\partial v^{\widetilde{j}}}{\partial x_i} \xi(|y|)V(x,y)dx\wedge dy+D\left(R_1\right)  \quad \text { for } R>R_2 ,
    \end{equation*}
    then
    $$
    Z'(R,\delta)=\int_{C^{(1)}_{R,\delta}} \bar{h}_{{\widetilde{i}} {\widetilde{j}}}(v) \frac{\partial v^{\widetilde{i}}}{\partial x_i} \frac{\partial v^{\widetilde{j}}}{\partial x_i} \xi(|y|)V(x,y)dS_{g_{can}}  ,
    $$
    where the derivative is taken with respect to $R$.
    By the Schwarz inequality, we have
    \begin{equation}\label{eq54}
        \begin{aligned}
            & \int_{C^{(1)}_{R,\delta}} \bar{h}_{{\widetilde{i}} {\widetilde{j}}}(v) \frac{\partial v^{\widetilde{i}}}{\partial x_i}  v^{\widetilde{j}} \nu^i  \xi(|y|)V(x,y)dS_{g_{can}}  \\
             \leq  & C \sqrt{\left(\int_{C^{(1)}_{R,\delta}} \bar{h}_{{\widetilde{i}} {\widetilde{j}}}(v) \frac{\partial v^{\widetilde{i}}}{\partial x_i} \frac{\partial v^{\widetilde{j}}}{\partial x_i}  \xi(|y|)V(x,y)dS_{g_{can}}\right)} \\
             & \times  \sqrt{\left(\int_{C^{(1)}_{R,\delta}} \bar{h}_{{\widetilde{i}} {\widetilde{j}}}(v) v^{\widetilde{i}} v^{\widetilde{j}}  \xi(|y|)V(x,y)dS_{g_{can}}\right)} .
        \end{aligned}
    \end{equation}
    Here, $C $ denotes an absolute constant that may vary with context.

    Furthermore, since the horizontal energy $E_{\mathcal{H}, \widetilde{\mathcal{H}}}(u)$ is infinite (by Lemma \ref{lm:mixdd} and $df_{\mathcal{H}, \widetilde{\mathcal{H}}}\neq 0$), there exists $R_3$ and $\delta_0$ such that $Z(R,\delta) > 0$ for $R \geq R_3$ and $\delta>\delta_0/2$. Then $\eqref{eq71}$ and $\eqref{eq53}$ imply
    \begin{equation*}
        Z(R,\delta)\leq \int_{C^{(1)}_{R,\delta}} 2 \bar{h}_{{\widetilde{i}} {\widetilde{j}}}(v) \frac{\partial v^{\widetilde{i}}}{\partial x_i}  v^{\widetilde{j}}   \nu^i\xi(|y|)V(x,y)dS_{g_{can}}.
    \end{equation*}
    Combining with $\eqref{eq54}$, we have 
    \[
    Z(R,\delta)^2 \leq CZ'(R,\delta) \left\{ \int_{C^{(1)}_{R,\delta}} \bar{h}_{{\widetilde{i}}{\widetilde{j}}}(v)v^{\widetilde{i}} v^{\widetilde{j}}  \, \xi(|y|)V(x,y)dS_{g_{can}} \right\} \quad \text{ for } R > R_3 \; \text{and}\; \delta>\delta_0/2.
    \]
    Denote
    \[
    M(R,\delta ) = \int_{C^{(1)}_{R,\delta}} \bar{h}_{{\widetilde{i}}{\widetilde{j}}}(v) v^{\widetilde{i}} v^{\widetilde{j}}  \xi(|y|) \, V(x,y)dS_{g_{can}}.
    \]
    Then for $R_4 \geq R \geq R_3$,
    \[
    \int_R^{R_4} \left(-\frac{1}{Z(r,\delta)}\right)' \, dr \geq C \int_R^{R_4} \frac{1}{M(r,\delta )} \, dr.
    \]
    
    On the other hand, by the assumption $( A_2)$ and the fact that $C_0(y)$ is bounded on $C^{(1)}_{R,\delta} $ , we have
    \begin{equation}\label{eq125}
        M(R,\delta )\leq C  \lambda (R)\int_{C^{(1)}_{R,\delta}}   V(x,y) \, dS_{g_{can}} 
        \leq C \lambda(R) R^{\sigma +1}\delta^k,
    \end{equation}
    where $\lambda$ can be chosen as a non-increasing function on $(R_3,\infty )$ such that  $\lambda(r)\geq \max_{|x|=r} \bar{h}_{{\widetilde{i}}{\widetilde{j}}}(v) v^{\widetilde{i}} v^{\widetilde{j}}$,     and $\lambda(r)\to 0$ as $r\to \infty$.
    It follows that 
    \begin{equation*}
        \frac{1}{Z(R,\delta)}\geq \frac{1}{Z(R,\delta)}-\frac{1}{Z(R_4,\delta)} \geq  C \int_R^{\infty} \frac{1}{M(r,\delta )} \, dr\geq C (\lambda (R) R^{\sigma }\delta^k)^{-1} ,
    \end{equation*}
    and thus
    \begin{equation}\label{eq:30}
        Z(R,\delta )\leq C\lambda(R)R^\sigma \delta^k  \quad\mathrm{for}\quad R\geq R_3 .
    \end{equation}
    Consequently,
    \begin{equation*}
        \begin{aligned}
            \int_{D_{R,\delta/2 }} e_{\mathcal{H}, \widetilde{\mathcal{H}}}(u)dV_g  
			&\leq \int_{D_{R, \delta} } e_{\HH}(u) \xi(|y|)dV_g \\
            &\leq C \lambda(R) R^{\sigma }\delta^k- D(R_1)+\int_{D_{R_2, \delta} }  e_{\HH}(u) \xi(|y|)dV_g\\
            &=C\Bigg(\lambda(R)\delta ^k+\frac{c(u)}{R^\sigma}\Bigg)R^\sigma,
        \end{aligned}
    \end{equation*}
where $c(u)$ is a constant depending on  $u$.
    
    Combining with Lemma \ref{lm:mixdd}, we deduce that
    \begin{equation}\label{eq:33}
        \begin{aligned}
            C\Bigg(\lambda(R)\delta ^k+\frac{c(u)}{R^\sigma}\Bigg) \geq & R^{-\sigma}\int_{D_{R,\delta}}e_{\mathcal{H},\widetilde{\mathcal{H}}}(u)dV_g\\
            \geq & \frac{1}{\rho^{\sigma }} \int_{D_{\rho,\delta}}e_{\mathcal{H},\widetilde{\mathcal{H}}}(u)dV_g
        \end{aligned} 
    \end{equation}
    for $R$ sufficiently large. Letting $R \to \infty $ and choosing $\rho$ large enough, we obtain
    \begin{equation*}
        0 \geq \frac{1}{\rho^{\sigma }} \int_{D_{\rho,\delta}}e_{\mathcal{H},\widetilde{\mathcal{H}}}(u)dV_g >0,
    \end{equation*}
    which is a contradiction.
\end{proof}
\begin{remark}
      If assumptions $( A_1)$  and  $( A_2)$ hold for different constants $\sigma$ and $\sigma'$,    then $\eqref{eq:33}$ becomes
        \begin{equation*}
            \begin{aligned}
                C\Bigg(\lambda(R)\delta ^k+\frac{c(u)}{R^{\sigma'}}\Bigg) \geq & R^{-\sigma '}\int_{D_{R,\delta}}e_{\mathcal{H},\widetilde{\mathcal{H}}}(u)dV_g\\
                \geq & \frac{R^{\sigma-\sigma '}}{\rho^{\sigma }} \int_{D_{\rho,\delta}}e_{\mathcal{H},\widetilde{\mathcal{H}}}(u)dV_g.
            \end{aligned} 
        \end{equation*}
        Thus the conclusion still holds provided $\sigma \geq \sigma '$.
\end{remark}
\begin{remark}
    Suppose $V$ is independent of $x$. Combining with assumption $(A_1)$,   $\eqref{eq125}$ becomes
          \begin{equation*}
              M(R,\delta )\leq C \lambda (R)\int_{C^{(1)}_{R,\delta}}   V(y) \, dS_{g_{can}} \leq C \Lambda(R,\delta) R^{m-1},
          \end{equation*}
          for some function $\Lambda(R,\delta).$  Therefore, the conclusion  remains valid provided $\sigma \geq m-2 $. 
  \end{remark} 

\begin{corollary}
    Let $u : (\mathbb{R}^{m+k}, \mathbb{R}^k, g) \rightarrow (N, \widetilde{\mathcal{F}}, \widetilde{g})$ be a horizontally harmonic map with $m > 2, k\geq 1$ and $u(x,y)\to q \in N$ as $(x,y)\to \infty$.
     If $g=C_0 g_{can}^{h}+ \eta(y)g_{can}^v$ for some constant $C_0>0$ and $\eta (y)$, then $du_{\mathcal{H}, \widetilde{\mathcal{H}}}= 0$.
\end{corollary}
\begin{proof}
    Since  assumption $( A_1)$ holds for $\sigma=m-2$ in this case, the claim follows from the above remark. 
\end{proof}

\begin{remark}
Note that the equation in the assumption $(A_1) $ is equivalent to 
$$  \frac{\partial V(x,y)}{\partial x^i}x^i \geq (2-m+\sigma)V(x,y). $$ 
\end{remark}

\begin{example}
    Both $( A_1)$ and  $( A_2)$ are satisfied in the following cases: 
    \begin{enumerate}[(i)]
        \item  $\phi =C|x|^{\frac{2(2-m+\sigma ) }{m-2}}$ and $\eta(x,y) =\eta (y)$ for some constants $C>0$ and $\sigma \geq m-2$.
        \item $\phi =C_1$  and $\eta(x,y) =C_2(y)|x|^{\frac{2(\sigma -m+2)}{k}}$ for some constants  $C_1>0$, $\sigma \geq m-2$ and some positive smooth function $C_2(y)$. This  also includes some warped products as a special case.  
    \end{enumerate}
    In particular, if $\phi =\eta \equiv  1 $, we get a Jin-type result for horizontally harmonic maps from the standard Euclidean spaces.
\end{example}
~

\subsection{Transversally harmonic maps}
~

Following the discussions in \S \ref{transverHarm} and \S \ref{se:Mixedconformal}, we may prove the following theorem in a similar way to Theorem \ref{thjin}.
\begin{theorem}
	Let $M^{m+k}$ be a complete Riemannian manifold with a pole $x_0$ and $r(x) $ be the distance from $x_0$.  
	 Suppose the radial curvature $K_r$ of $M$ satisfies  $-\alpha ^2\leq  K_r \leq -\beta ^2  $, where $\alpha ,\beta >0 $ and $(m+k-1)\beta -2\alpha \geq 0$. 
     Let $u:\foli \to (N^{n+l}, \widetilde{\mathcal{F}}^l , \widetilde{g})$ be a transversally harmonic map with $m>2,k \geq 1$ and  
     $$ \dist_N (u(x),q) \leq \eta (r(x)) \quad \text{for} \quad r \gg 1, $$ 
     for some point $q\in N$, where $\eta $ is a non-increasing function with  
	$$ \eta(t)=O(t^{\frac{\lambda }{2}+\epsilon '} e^{-(m+k-1)\frac{\alpha }{2} t}) \quad \text{as} \quad t \to \infty $$ 
	for some $\frac{1}{2}>\epsilon ' >0$ and $\lambda =m+k-\frac{2 \alpha }{\beta }.$  
	Then $u$ is horizontally constant. 
\end{theorem}

\begin{proof}
    Let $(U,\varphi ,\widetilde{x}^{\widetilde{A}})$, $(\widetilde{U},
\widetilde{\varphi },\widetilde{x}^{\widetilde{i}})$ be coordinate charts, with $\pi :U\rightarrow 
\widetilde{U}$ and  $\overline{q}=\pi (q)$  as defined in \S \ref{se:Mixedconformal}. Since $%
u(x)\rightarrow q$ as $x\rightarrow \infty $, there is an $R_{1}>0$ such that if $r(x)>R_{1}$, then $u(x)\in U$. Define the map $v:M\backslash
B_{R_{1}}\rightarrow U$ by setting $v:=\pi \circ u$. 

Let $\{e_{A}\}_{A=1}^{m+k}$ be an adapted frame field around any point $p\in
M\backslash B_{R_{1}}$. For any $R>R_{1}$, the horizontal energy of $u$ on $%
B_{R}$, denoted by $E_{\mathcal{H},\widetilde{\mathcal{H}}}^{R}(u)$, is given by  
	$$ \begin{aligned}
        E_{\mathcal{H},\widetilde{\mathcal{H}}}^{R}(u) =& \int_{B_{R_1}} e_{\HH}(u)dV_g+ \frac{1}{2}\int_{B_R\setminus B_{R_1}}  h( du_{\HH}(e_k),du_{\HH}(e_k) )  dV_g\\
	=& \int_{B_{R_1}} e_{\HH}(u)dV_g+\frac{1}{2}\int_{B_R \setminus B_{R_1}} \bar h(  d v_{\mathcal{H} }(e_k), d v_{\mathcal{H}}(e_k))  dV_g ,
\end{aligned}  $$
where $\bar{h} $ is the metric on $ \widetilde{U}$.  
 For simplicity, set
	$$ E^R_{\mathcal{H}, \cdot  }(v)=: \frac{1}{2}\int_{B_R \setminus B_{R_1}} \bar h(  d v_{\mathcal{H} }(e_k), d v_{\mathcal{H} }(e_k))  dV_g \quad \text{for } R >R_1.  $$
 Writing   $d v_{\mathcal{H}}(e_k)=v^{\widetilde{i}}_k \frac{\partial }{\partial x^{\widetilde{i}}}  $ and $\bar{h}_{\widetilde{i}\widetilde{j}}= \bar{h}(\frac{\partial }{\partial x^{\widetilde{i}}},\frac{\partial }{\partial x^{\widetilde{j}}}  ) $, we have   
$$ E^R_{\mathcal{H}, \cdot }(v)=\frac{1}{2}\int_{B_R \setminus B_{R_1} }  v^{\widetilde{i}}_k v^{\widetilde{j}}_k \bar h_{\widetilde{i}{\widetilde{j}}}(v)  dV_g. $$
As in Theorem \ref{thjin}, by choosing a sufficiently large $R_1$, we get  
	\begin{equation}\label{eq4900}
		\left( \frac{\partial \bar h_{{\widetilde{i}}{\widetilde{j}}}(v)}{\partial x^{\widetilde{k}}} v^{\widetilde{k}} + 2h_{{\widetilde{i}}{\widetilde{j}}}(v) \right) \geq (h_{{\widetilde{i}}{\widetilde{j}}}(v)) \quad \text{ for } r(x) > R_1
	\end{equation}
in sense of matrices. 
For any $\widetilde{w}  \in C^{2}_0(M\backslash B_{R_1},\varphi (U))$ and sufficiently small $t$ , we  consider the variation $u+t\widetilde{w}:M \to N$ defined as follows:
$$ (u+t \widetilde{w})(x)= \begin{cases}
	 u(x)& \text{if} \quad x \in B_{R_1}\\
	 \varphi ^{-1}[\varphi \left( u(x)\right)+t\widetilde{w}(x)]&\text{if} \quad x \in M\backslash B_{R_1}.
\end{cases}  $$  
By the definition of transversally harmonic maps,  we have
	\[
	\left. \frac{d}{dt} \right|_{t=0} E_{\mathcal{H}, \widetilde{\mathcal{H}}}(u + t \widetilde{w}) = 0,
	\]
	that is, by setting $w=\Pr_n \circ \widetilde{w}\in C_{0}^{2}(M \setminus B_{R_1}, \widetilde{\varphi } (\widetilde{U})) $,
	$$ \begin{aligned}
		0=&\int_{M\backslash B_{R_1}}  \left( 2 \bar h_{{\widetilde{i}}{\widetilde{j}}}(v)v^{\widetilde{i}}_k w^{\widetilde{j}}_k + \frac{\partial \bar h_{{\widetilde{i}}{\widetilde{j}}}(v)}{\partial x^{\widetilde{l}}} w^{\widetilde{l}} v^{\widetilde{i}}_k v^{\widetilde{j}}_k \right) dV_g .
	\end{aligned}  $$
	Choose \( \widetilde{w}(x) = \psi(r(x)) \varphi (u(x)) \)  for  \( \psi(t) \in C_0^\infty(R_1, \infty) \). We have
	\begin{equation}\label{eq500}
		\begin{aligned}
			&\int_{M\backslash B_{R_1}} \left( 2\bar{h}_{{\widetilde{i}}{\widetilde{j}}}(v) + \frac{\partial \bar{h}_{{\widetilde{i}}{\widetilde{j}}}(v)}{\partial x^{\widetilde{l}}} v^{\widetilde{l}} \right)  v^{\widetilde{i}}_k v^{\widetilde{j}}_k  \psi dV_g\\
			=& - \int_{M \backslash B_{R_1}} 2\bar{h}_{{\widetilde{i}}{\widetilde{j}}}(v)v^{\widetilde{i}}_k  v^{\widetilde{j}} \psi_{k}  dV_g .
		\end{aligned}
	\end{equation}
	As before, for $0<\varepsilon \leq 1$, we define $\zeta_{\varepsilon}$ as follows:
	$$
	\zeta_{\varepsilon}(t)= \begin{cases}1 & t \leq 1; \\ 1+\frac{1-t}{\varepsilon} & 1<t<1+\varepsilon ; \\ 0 & t \geq 1+\varepsilon .\end{cases}
	$$
	In $\eqref{eq500}$, choose
	$$
	\psi(r(x))=\zeta_{\varepsilon}\left(\frac{r(x)}{R}\right)\left(1-\zeta_1\left(\frac{r(x)}{R_1}\right)\right) 
	$$
	for some $R> R_2 $,  where $R_2:=2R_1$.
	By noting that
	$$
	e_k(\zeta_{\varepsilon}\left(\frac{r(x)}{R}\right))=-\frac{1}{R \varepsilon} r_k  \quad  \text { for } R<r(x)<R(1+\varepsilon),
	$$
	we derive 
	\begin{equation}\label{eq520}
		\begin{aligned}
			& \int_{B_{R}\backslash B_{R_2}}\left(2 \bar{h}_{{\widetilde{i}} {\widetilde{j}}}(v)+\frac{\partial \bar{h}_{{\widetilde{i}} {\widetilde{j}}}(v)}{\partial x^{\widetilde{l}}} v^{\widetilde{l}}\right)  v^{\widetilde{i}}_k v^{\widetilde{j}}_k dV_g+D\left(R_1\right) \\
			&+\int_{M\backslash B_{R}} \left( 2\bar{h}_{{\widetilde{i}}{\widetilde{j}}}(v) + \frac{\partial \bar{h}_{{\widetilde{i}}{\widetilde{j}}}(v)}{\partial x^{\widetilde{l}}} v^{\widetilde{l}} \right)  v^{\widetilde{i}}_k v^{\widetilde{j}}_k \zeta_{\varepsilon}\left(\frac{r(x)}{R}\right)dV_g\\
			=&\frac{1}{R \varepsilon}\int_{B_{(1+\epsilon )R}\backslash B_{R}} 2 \bar{h}_{{\widetilde{i}} {\widetilde{j}}}(v) v^{\widetilde{j}} v^{\widetilde{i}}_k r_k   dV_g,
		\end{aligned}
	\end{equation}
	where
	$$
	\begin{aligned}
		D\left(R_1\right)= & \int_{B_{R_2}\backslash B_{R_1}}\left(2 \bar{h}_{{\widetilde{i}} {\widetilde{j}}}(v)+\frac{\partial \bar{h}_{{\widetilde{i}} {\widetilde{j}}}(v)}{\partial x^{\widetilde{l}}} v^{\widetilde{l}}\right)  v^{\widetilde{i}}_k v^{\widetilde{j}}_k \left(1-\zeta_1\left(\frac{r(x)}{R_1}\right)\right)  dV_g\\
		& -\frac{1}{R_1}\int_{B_{(1+\varepsilon )R_1}\backslash B_{R_1}} 2 \bar{h}_{{\widetilde{i}} {\widetilde{j}}}(v) v^{\widetilde{j}}v^{\widetilde{i}}_k  r_k(x)  dV_g.
	\end{aligned}
	$$
	By sending $\varepsilon \to 0$ in $\eqref{eq520}$, we have 
	\begin{equation}\label{eq530}
		\begin{aligned}
			& \int_{B_{R}\backslash B_{R_2}}\left(2 \bar{h}_{{\widetilde{i}} {\widetilde{j}}}(v)+\frac{\partial \bar{h}_{{\widetilde{i}} {\widetilde{j}}}(v)}{\partial x^{\widetilde{l}}} v^{\widetilde{l}}\right) v^{\widetilde{i}}_k v^{\widetilde{j}}_k  dV_g+D\left(R_1\right) \\
			=& \int_{\partial B_R} 2 \bar{h}_{{\widetilde{i}} {\widetilde{j}}}(v) v^{\widetilde{j}} v^{\widetilde{i}}_k r_k    dS_g.
		\end{aligned}
	\end{equation}
	Set
	\begin{equation*}
		Z(R)=\int_{B_{R}\backslash B_{R_2}} \bar{h}_{{\widetilde{i}} {\widetilde{j}}}(v) v^{\widetilde{i}}_k v^{\widetilde{j}}_k   dV_g+D\left(R_1\right)  \quad \text { for } R>R_2 .
	\end{equation*}
	Then
	$$
	Z'(R)=\int_{\partial B_R} \bar{h}_{{\widetilde{i}} {\widetilde{j}}}(v) v^{\widetilde{i}}_k v^{\widetilde{j}}_k    dS_g  .
	$$
	By the Schwarz inequality, we have
	\begin{equation}\label{eq540}
		\begin{aligned}
			 \int_{\partial B_R} 2 \bar{h}_{{\widetilde{i}} {\widetilde{j}}}(v) v^{\widetilde{j}} v^{\widetilde{i}}_k r_k     dS_g  \leq & \int_{\partial B_R} |2 \bar{h}_{{\widetilde{i}} {\widetilde{j}}}(v) v^{\widetilde{j}} v^{\widetilde{i}}_k |     dS_g  \\
			 \leq & C \sqrt{\left(\int_{\partial B_R} \bar{h}_{{\widetilde{i}} {\widetilde{j}}}(v) v^{\widetilde{i}}_k v^{\widetilde{j}}_k  dS_g\right)} \times \sqrt{\left(\int_{\partial B_R} \bar{h}_{{\widetilde{i}} {\widetilde{j}}}(v) v^{\widetilde{i}} v^{\widetilde{j}}  dS_g\right)} \\
			 = &C \sqrt{ Z'(R)} \times \sqrt{\left(\int_{\partial B_R} \bar{h}_{{\widetilde{i}} {\widetilde{j}}}(v) v^{\widetilde{i}} v^{\widetilde{j}}  dS_g\right)}.
		\end{aligned}
	\end{equation}

	Moreover, if $du_{\mathcal{H}, \widetilde{\mathcal{H}}}\neq 0$, the horizontal energy $E_{\mathcal{H}, \widetilde{\mathcal{H}}}(u)$ is infinite by Lemma \ref{lm:gwtype}. Thus, there exists an $R_3$, such that $Z(R) > 0$ for $R \geq R_3$. Then $\eqref{eq4900}$,$\eqref{eq530}$ 
    and $\eqref{eq540}$ imply 
		\[
		Z(R)^2 \leq CZ'(R) \left\{ \int_{\partial B_R} \bar{h}_{{\widetilde{i}}{\widetilde{j}}}(v)v^{\widetilde{i}} v^{\widetilde{j}}  \, dS_g \right\} \quad \text{ for } R > R_3.
		\]
		Set
		\[
		M(R) = \int_{\partial B_R} \bar{h}_{{\widetilde{i}}{\widetilde{j}}}(v) v^{\widetilde{i}} v^{\widetilde{j}}  dS_g.
		\]
		For $R_4 \geq R \geq R_3$, we obtain
		\[
		\int_R^{R_4} \left(-\frac{1}{Z(r )}\right)' \, dr \geq C \int_R^{R_4} \frac{1}{M(r)} \, dr.
		\]
	
		On the other hand, by the definition of $\eta $ and the volume comparison theorem, we have
		\begin{equation}\label{eq1250}
			M(R)\leq C |\partial B_R| \eta^2 (R) \leq C e^{(m+k-1)\alpha R} \eta^2 (R).
		\end{equation}
		It follows that 
		\begin{equation*}
			\frac{1}{Z(R )}\geq \frac{1}{Z(R )}-\frac{1}{Z(R_4 )} \geq  C \int_R^{\infty} \frac{1}{M(r)} \, dr\geq C \int_R^{\infty} (e^{(m+k-1)\alpha r} \eta^2 (r))^{-1}dr .
		\end{equation*}
	Since $\lambda \geq 1$, we have 
	   $$ \int_R^{\infty} (e^{(m+k-1)\alpha r} \eta^2 (r))^{-1}dr \geq C \int_R^{\infty} (t^{\lambda +2 \epsilon '})^{-1}dr= C(\lambda ,\epsilon ') R^{-\lambda -2 \epsilon '+1}. $$ 
	Thus
		\begin{equation}\label{eq:300}
			Z(R  )\leq  C(\lambda ,\epsilon ') R^{\lambda +2 \epsilon '-1} .
		\end{equation}
	Clearly,
		\begin{equation*}
			\begin{aligned}
				E^R_{\HH}(u)  &=\frac{1}{2}Z(R)-\frac{1}{2}D(R_{1})+\int_{B_{R_2}} e_{\HH}(u)dV_g\\
				&\leq  \frac{1}{2}C(\lambda ,\epsilon ') R^{\lambda +2 \epsilon ' -1} -\frac{1}{2}D(R_1)+\int_{B_{R_2} }e_{\HH}(u)dV_g\\
                &=C\left(R^{2 \epsilon '-1}+\frac{c(u)}{R^\lambda }\right)R^\lambda ,
			\end{aligned}
		\end{equation*}
		where $c(u)$ is a constant depending on  $u$.
		
		Combining with Lemma \ref{lm:gwtype} and $\epsilon '<\frac{1}{2}$, we  have
		\begin{equation}
			\begin{aligned}
				C\left(R^{2 \epsilon '-1}+\frac{c(u)}{R^\lambda }\right)  \geq & R^{- \lambda}\int_{B_{R}} e_{\HH}(u) dV_g\\
				\geq & \frac{1}{\rho^{ \lambda }} \int_{B_\rho }e_{\HH}(u)dV_g
			\end{aligned} 
		\end{equation}
		for sufficiently large $R$. By fixing a large $\rho $  and	letting $R \to \infty $, we obtain
		\begin{equation*}
			0 \geq \frac{1}{\rho^{ \lambda }} \int_{B_{\rho}}e_{\HH}(u)dV_g >0,
		\end{equation*}
		which is a contradiction.
	\end{proof}
\begin{remark}
    This theorem establishes a Jin-type result under the assumption of a specific convergence rate towards the point $q$,
 a method originally introduced by \cite{rs2000energy}. 
\end{remark}

\section*{Acknowledgments}
Tian Chong was supported by National Natural Science Foundation of China (NSFC) No. 12001360. 
Yuxin Dong was supported by NSFC Grant No. 12171091.
Xin Huang was supported by the Startup Foundation for Introducing Talent of NUIST.   
and Hui Liu was supported by the China Scholarship Council (No. 202306100156).
The authors would like to thank Prof. V. Branding for his valuable comments.

\end{document}